\documentclass[12pt]{amsart}
\usepackage{amssymb,latexsym,eufrak,amsmath,amscd, graphicx}
\usepackage[dvipsnames,usenames]{color}

\usepackage{amsfonts}
\usepackage[all]{xypic}
\UseTips

\renewcommand{\to}[1][]{\xrightarrow{\ #1\ }}

\newcommand{\forget}[1]{}  

\makeatletter
\renewcommand{\theenumi}{\@roman\c@enumi}
\makeatother

\renewcommand{\phi}{\varphi}
\renewcommand{\epsilon}{\varepsilon}
\renewcommand{\theta}{\vartheta}

\newcommand{\llbracket}{[\negthinspace[}
\newcommand{\rrbracket}{]\negthinspace]}

\def\ZZ{{\mathbf Z}}

\def\AAA{{\mathbf A}}

\def\QQ{{\mathbf Q}}

\def\cO{\mathcal{O}}

\def\fra{\mathfrak{a}}
\def\frb{\mathfrak{b}}
\def\frc{\mathfrak{c}}

\def\frm{\mathfrak{m}}
\def\frq{\mathfrak{q}}
\def\frp{\mathfrak{p}}

\def\Id{{\mathcal Id}_{\bf Q}} 

\def\o{\circ}

\def\.{\cdot}
\def\({\Big{(}}
\def\){\Big{)}}
\def\^{\widehat}
\def\~{\widetilde}
\def\*{{}^*\!}
\def\[{\llbracket}
\def\]{\rrbracket}

\renewcommand{\and}{ \quad \text{and} \quad }

\DeclareMathOperator{\codim}{codim}

 \DeclareMathOperator{\Spec}{Spec}
 
 \DeclareMathOperator{\lct}{lct}
 
 \DeclareMathOperator{\ord}{ord}

\newtheorem{lemma}{Lemma}[section]
\newtheorem{theorem}[lemma]{Theorem}
\newtheorem{corollary}[lemma]{Corollary}
\newtheorem{proposition}[lemma]{Proposition}

\newtheorem{conjecture}[lemma]{Conjecture}

\theoremstyle{definition}
\newtheorem{definition}[lemma]{Definition}
\newtheorem{remark}[lemma]{Remark}

\newtheorem{example}[lemma]{Example}
\newtheorem{question}[lemma]{Question}

\theoremstyle{remark}
\newtheorem*{remark*}{Remark}
\newtheorem*{note*}{Note}

\begin{document}

\title{Toward an inductive description of singularities of pairs}

\author[M. Musta\c{t}\u{a}]{Mircea~Musta\c{t}\u{a}}
\address{Department of Mathematics, University of Michigan,
Ann Arbor, MI 48109, USA} \email{{\tt mmustata@umich.edu}}

\markboth{M.~Musta\c t\u a}{Inductive description of singularities
of pairs}

\begin{abstract}
Motivated by Shokurov's ACC Conjecture for log canonical thresholds, we propose 
an inductive point of view on singularities of pairs, in the case when the ambient variety
is smooth.
Our main result characterizes the log canonicity of a pair in dimension $n+1$ by the log canonicity of another pair (not effective, in general) in dimension $n$. 
\end{abstract}

\thanks{2000\,\emph{Mathematics Subject Classification}.
 Primary 14B05; Secondary 14B10, 14E30.
\newline The author was partially supported by
 NSF grant DMS-0758454 and
  a Packard Fellowship}

\maketitle

\section{Introduction}

Singularities in higher-dimensional birational geometry are classified according to orders
of vanishing along divisorial valuations. A key tool in the study of singularities is provided
by log resolutions: one can show that instead of dealing with all possible divisorial valuations, 
it is enough to consider the ones that correspond to divisors on a given log resolution.
Our goal in this note is to show that one can try to approach singularities by induction on the dimension
of the ambient space. The main ingredient for doing this is the characterization of singularities of pairs
in terms of spaces of arcs, as developed in \cite{Mus} and \cite{ELM}.

We work in the following setup. Let $K$ be an algebraically closed field of characteristic zero.
Our ambient scheme $X$ will always be nonsingular: it will either be a smooth variety over $K$,
or, most of the time, ${\rm Spec}\,K\[x_1,\ldots,x_n\]$.  The usual setting for singularities in birational geometry is that of pairs $(X,D)$, where $D$ is a divisor on $X$ with rational coefficients. For us it is more natural to allow arbitrary ideals, hence our pairs will be of the form $(X,J_1^{q_1}J_2^{-q_2})$,
where $J_1$ and $J_2$ are nonzero ideals on $X$, and $q_1$ and $q_2$ are nonnegative rational numbers. We say that such a pair is effective if $q_2=0$.
The condition for a pair $(X,J_1^{q_1}J_2^{-q_2})$ to be \emph{log canonical} is typically given
in terms of the orders of vanishing of $J_1$ and $J_2$ along the divisorial valuations of the function field of $X$. 
If $J$ is a proper nonzero ideal on $X$, then the \emph{log canonical threshold} $\lct(J)$ is the largest 
$q>0$ such that the pair $(X,J^q)$ is log canonical. This invariant plays an important role in birational geometry. For other points of view on the log canonical threshold, and for some applications, see \cite{EM}.

Inductive arguments have always been present in the theory of singularities of pairs 
through the adjunction formula, which leads to Adjunction and Inversion of Adjunction type results
(see \S 7 in \cite{kollar}). However, these results apply to only special pairs, while we are interested
in general ones. 
By an inductive description of singularities of \emph{effective} pairs we mean a positive answer to the following question:

\begin{question}\label{inductive_description1}
Given integers $n, d\geq 1$ and a positive rational number $c$, consider formal power series of the form
$$f=y^d+\sum_{i=1}^da_i(x_1,\ldots,x_n)y^{d-i},$$
with $a_i\in K\[x_1,\ldots,x_n\]$, of positive order.
Can one find finitely many polynomials $P_i\in K[z_1,\ldots,z_d]$, with $1\leq i\leq r$, and a positive 
rational number $q$, such that for every
$f$ as above we have
$$\left({\rm Spec}\,K\[x_1,\ldots,x_d,y\],f^c\right)\,\text{is log canonical if and only if}$$
$$\left({\rm Spec}\,K\[x_1,\ldots,x_n\], (P_1(a_1,\ldots,a_d),\ldots,P_r(a_1,\ldots,a_d))^{q}\right)
\,\text{is log canonical}\,?$$
\end{question}

Our main motivation for looking at this question comes from Shokurov's ACC Conjecture on log canonical thresholds. This is one of the most intriguing open problems in the field.

\begin{conjecture}\label{ACC}${\rm (}$\cite{Sho}${\rm )}$
If $n\geq 1$ is fixed, the set
$${\mathcal HT}_n:=\{\lct(X,D)\mid X\,\text{nonsingular},\,\dim(X)=n,\,D\,\text{effective divisor}\}$$
satisfies ACC, that is, it contains no strictly increasing sequences.
\end{conjecture}

Several proofs of the conjecture are known for $n=2$. Shokurov proved it in \cite{Sho} using 
the Minimal Model Program. Phong and Sturm gave an analytic proof in \cite{PS}, and
Favre and Jonsson proved it in \cite{FJ} using their valuation tree. The case $n=3$ was proved by
Alexeev in \cite{alexeev}. We should point out that one can consider singularities of pairs when the
ambient variety is not assumed to be nonsingular, and it is in this setting that the above conjecture was made. In this stronger form, it is related to Termination of Flips, the major open problem 
left in the Minimal Model Program (see \cite{birkar} for the precise statement). 

It was shown in \cite{dFM} that Conjecture~\ref{ACC} is equivalent to the following statement:
given a positive rational number $q$, there is $N$ such that for every 
$f, g\in K\[x_1,\ldots,x_n\]$ with $\ord(f)\geq 1$ and $\ord(g)\geq N$,
we have $\lct(f)\geq q$ if and only if
$\lct(f+g)\geq q$. Using the same ideas, we show in Proposition~\ref{implication}  below
that a positive answer to Question~\ref{inductive_description1} in dimension $n$, plus Conjecture~\ref{ACC} in dimension $n$ would imply the conjecture
in dimension $n+1$. The idea is that in an increasing sequence of log canonical thresholds, 
the orders of the corresponding equations are bounded above, hence we may assume these orders to be constant.
In this case, by
Weierstarss preparation we can write the equations as polynomials of constant degree in the last variable, in which case the ideas in \cite{dFM}, and the positive 
answer to Question~\ref{inductive_description1} allow us to apply induction.

Our main result is in the direction of Question~\ref{inductive_description1}, though we can not
remain in the setting of effective pairs.

\begin{theorem}\label{main_intro}
Fix an integer $d\geq 1$ and a positive rational number $c$, and consider formal power series of the form
$$f=y^d+\sum_{i=1}^da_i(x_1,\ldots,x_n)y^{d-i},$$
with $n\geq 1$, and $a_i\in K\[x_1,\ldots,x_n\]$ of positive order.  There are finitely many polynomials $P_i, Q_j\in \QQ[z_1,\ldots,z_d]$, with $1\leq i\leq r$ and $1\leq j\leq s$, and nonnegative 
rational numbers $p$ and $q$ ${\rm (}$everything depending only on $d$ and $c$${\rm )}$, such that for every
$f$ as above we have 
$$\left({\rm Spec}\,K\[x_1,\ldots,x_n,y\],f^c\right)\,\text{is log canonical if and only if}$$
$$\left({\rm Spec}\,K\[x_1,\ldots,x_n\], I^{p}J^{-q}\right)
\,\text{is log canonical},$$
where $I=(P_1(a_1,\ldots,a_d),\ldots, P_r(a_1,\ldots,a_d))$ and 
$J=(Q_1(a_1,\ldots,a_d),\ldots,Q_s(a_1,\ldots,a_d))$. 
\end{theorem}

In fact, we will obtain the $P_i$ and the $Q_j$ in a rather explicit manner (see Theorem~\ref{main}
below for the precise statement). Let us consider the special case $d=3$ (see Example~\ref{degree_three} below).
After a linear change of variable, we may assume that $f=y^3+a(x)y+b(x)$, with 
$a,b\in K\[x_1,\ldots,x_n\]$ of positive order.
If $\Delta_f=4a(x)^3+27b(x)^2$ is the discriminant of $f$, then for every 
$c$ with $\frac{2}{3}\leq c\leq 1$, we have
$$ \lct(f)\geq c\,\,\text{if and only if}\,\,\left({\rm Spec}\,K\[x_1,\ldots,x_n\], (\Delta_f)^{c-\frac{1}{2}}
(a(x)^3,b(x)^2)^{-\frac{3c-2}{6}}\right)\,\,\text{is log canonical}.$$

While in Theorem~\ref{main_intro} the pair we obtain is not effective in general, 
the explicit description we obtain allows us to compare 
$I^{p}$ and $J^{q}$: we show that for every divisor $E$
over ${\rm Spec}\,K[z_1,\ldots,z_d]$ we have
$$p\cdot \ord_E(I)\geq \frac{d}{d-1}q\cdot\ord_E(J).$$
This is enough to yield a new proof of Conjecture~\ref{ACC} for $n=2$ (see Corollary~\ref{dimension_two} below). 

The fact that in Theorem~\ref{main_intro} we describe the log canonicity of the given effective pair in terms of the log canonicity of a possibly non-effective pair raises some doubts regarding an affirmative answer to Question~\ref{inductive_description1}. On the other hand, the theorem suggests that there are classes of non-effective pairs that behave like the effective ones (we recall that effective pairs have 
in general much better properties than non-effective ones do: for example, they satisfy semicontinuity of log canonical thresholds and Inversion of Adjunction). It is conceivable that applying to such an extended class of pairs our point of view might lead to a solution of Conjecture~\ref{ACC}.

Describing effective log canonical pairs in dimension $(n+1)$ in terms of possibly non-effective
pairs in dimension $n$ is not enough for a truly inductive approach to singularities of pairs.
 On the other hand, we believe that even if we start with a non-effective pair in dimension $(n+1)$, one should be able to characterize its log canonicity in terms of the log canonicity of other suitable pairs in dimension $n$ (though it might be hard to get such an explicit description as in Theorem~\ref{main} below). We hope to return to this aspect in future work.

\bigskip

The main tool in the proof of Theorem~\ref{main_intro} is the characterization of log canonical pairs in terms
of the codimensions of the contact loci associated to $f$. More precisely, given 
$f\in K\[x_1,\ldots,x_n,y\]$, the log canonical threshold of $f$ gives a lower bound for the codimension of the subsets
of $\left(tK\[t\]\right)^{n+1}$ in terms of the order of vanishing of $f$ along these subsets; see 
Theorem~\ref{description} below for the precise statement. Our idea is to first plug in
elements in $u\in (tK\[t\])^n$, and then  estimate the codimension of those $w$ such that
$\ord(f(u,w)=m$ by treating $f(u,y)$ as a polynomial with coefficients in $K\[t\]$, and considering
its factorization in the field of Puiseux formal power series. 

The paper is structured as follows. In the second section we give a brief review of the definition
of log canonical pairs in the setting that we need later in the paper. The following section is devoted to the characterization of log canonicity in terms of codimension of cylinders in the 
space of arcs. We here explain how to deduce this characterization in the case of formal power series
from the one in the case of smooth varieties, that is treated in \cite{ELM}. The fourth section treats
polynomials in $K\[t\][y]$. Given a monic such polynomial of degree $d$, we show how to estimate the orders of various expressions in its roots. This part is very elementary, but it provides the ingredients for our proof of Theorem~\ref{main_intro}. In the last section we explain the connection between
 Question~\ref{inductive_description1} and Conjecture~\ref{ACC}, and prove our main result
 on the inductive description of log canonical singularities.

\subsection{Acknowledgments}
We are grateful to Tommaso de Fernex, Lawrence Ein, and Rob Lazarsfeld for many 
discussions, and for their suggestions and insight on topics related to this paper.

\section{Review of log canonical pairs}

In this section we review the definition of log canonical pairs. Our main reference for this topic
is \cite{kollar}. However, our setting is slightly different from the one in \emph{loc. cit.}, since
we need to work with ideals in a formal power series ring. 
Extending most of the results to this setting poses no problem. Such an extension has been carried out
for the basic properties of the log canonical threshold in \cite{dFM}.

We start by recalling some terminology that will be useful for dealing with singularities of pairs.
In this approach we follow \cite{Kawakita}. Let $X$ be a connected Noetherian normal scheme. 
A \emph{divisor over $X$} 
is a prime divisor on some normal scheme $X'$, having a proper, birational morphism  $\phi\colon X'\to X$. Such a divisor
$E$ gives a valuation of the function field of $X$ that we denote by $\ord_E$. We identify two divisors
over $X$ if they give the same valuation. The center $c_X(E)$ of $E$ on $X$ is the image $\phi(E)$.

A $\QQ$-ideal
$\fra$ is an equivalence class of symbols $J^q$, where $J$ is a nonzero ideal on $X$,
and $q$ is a nonnegative rational number. The equivalence relation is given by
$J_1^{q_1}\sim J_2^{q_2}$ if and only if for some (every) positive integer $m$ with $mq_1, mq_2\in\ZZ$, the ideals
$J_1^{mq_1}$ and $J_2^{mq_2}$ have the same integral closure (we refer to \cite{positivity}, 
\S 9.6.A for basic facts about the integral closure of ideals). 
Equivalently, for every divisor $E$ over $X$, we have
$q_1\cdot\ord_E(J_1)=q_2\cdot\ord_E(J_2)$.  
It is convenient to 
consider also the zero ideal as a $\QQ$-ideal: by convention,
this is equal to $(0)^q$ for every nonnegative $q\in\QQ$.
We denote by $\Id(X)$ the set of $\QQ$-ideals on $X$.  

One can define the sum and the product of the $\QQ$-ideals
$J_1^{q_1},\ldots, J_r^{q_r}$,  as follows. If $m$ is a positive integer
such that   $mq_i\in\ZZ$ for every $i$, then
$$\prod_{i=1}^rJ_i^{q_i}:=\left(\prod_{i=1}^r J_i^{mq_i}\right)^{1/m},$$
$$\sum_{i=1}^rJ_i^{q_i}:=\left(\sum_{i=1}^rJ_i^{mq_i}\right)^{1/m}.$$
Furthermore, we can define rational powers of $\QQ$-ideals: if 
$\fra=J^q$, and if $q'$ is a positive rational number, 
then $\fra^{q'}=J^{qq'}$. 
It is straightforward to see that all these operations are well-defined.

The set of nonzero $\QQ$-ideals forms a semigroup with respect to the product of 
$\QQ$-ideals. This semigroup is integral, that is, $\fra_1\cdot\frb=\fra_2\cdot\frb$ implies $\fra_1=\fra_2$ in $\Id(X)$. Therefore 
it admits a universal semigroup morphism to an abelian group $\Id(X)^{\rm gp}$, and this morphism 
is injective. The elements in $\Id(X)^{\rm gp}$ are symbols $\fra\frb^{-1}$, with $\fra,\frb$ nonzero elements in
$\Id(X)$, and such that 
$\fra\frb^{-1}$ and $\fra_1\frb_1^{-1}$ give the same element in this group if and only if
$\fra\frb_1=\fra_1\frb$ in $\Id(X)$.

Given a divisor 
$E$ over $X$, and $\fra=J^q\in\Id(X)$, we define $\ord_E(\fra):=q\cdot\ord_E(J)$. The definition of
$\Id(X)$ was set up so that $\ord_E(\fra)$ is well-defined. Note that if $\fra$ is nonzero, then $\ord_E(\fra)$ is finite, and  $\ord_E$ extends to a group homomorphism $\Id(X)^{\rm gp}\to\QQ$.

If $f\colon Y\to X$ is a morphism of schemes as above, then we can pull-back
$\QQ$-ideals in the obvious way: $f^{-1}(J^{\lambda})=(J\cdot\cO_Y)^{\lambda}$.
Note that this pull-back can be the zero $\QQ$-ideal.
This operation is compatible with sums, products, and rational powers. 

We now introduce some notation that will be repeatedly used later.
Suppose that $K$ is a field of characteristic zero.
Given a $\QQ$-ideal $\frb$ in $\QQ[z_1,\ldots,z_d]$, and formal power series 
$a_1,\ldots,a_d\in K\[x_1,\ldots,x_n\]$,  we denote by
$\frb(a_1,\ldots,a_d)$ the pull-back of $\frb$ by the scheme morphism induced by
$\QQ[z_1,\ldots,z_d]\to K[x_1,\ldots,x_n]$, $z_i\to a_i$. 
Similarly, if $a_1,\ldots,a_d\in K\[t^{1/N}\]$, for some $N\geq 1$, then we 
get a $\QQ$-ideal $\frb(a_1,\ldots,a_n)$ in $K\[t^{1/N}\]$. Note that we can consider the order of this ideal:
the order of the $\QQ$-ideal $(t)^q$ in $K\[t^{1/N}\]$ is $q$ (we make the convention that 
$\ord(0)=\infty$).

More generally, if $\gamma\colon \Spec\,K\[t^{1/N}\]\to X$ is a scheme morphism, and if $\fra$
is a $\QQ$-ideal on $X$, then we denote by $\ord_{\gamma}(\fra)$ the order of 
$\gamma^{-1}(\fra)$. Note that we have
$$\ord_{\gamma}\left(\prod_{i=1}^r\fra_i\right)=\sum_{i=1}^r\ord_{\gamma}(\fra_i),\,\,
\ord_{\gamma}\left(\sum_{i=1}^r\fra_i\right)=\min_{i=1}^r\{\ord_{\gamma}(\fra_i)\},\,\,
\ord_{\gamma}(\fra^c)=c\cdot\ord_{\gamma}(\fra)$$
for every $\QQ$-ideals $\fra_i$ and $\fra$, and every positive rational number $c$.
Similar formulas hold if instead of $\ord_{\gamma}$ we consider $\ord_E$, for a divisor $E$
over $X$.

It is convenient to also consider an order relation on the nonzero elements of $\Id(X)$.
For two $\QQ$-ideals $\fra$ and $\frb$ on $X$ we put $\fra\leq \frb$ if 
$\ord_E(\fra)\geq\ord_E(\frb)$ for every divisor $E$ over $X$.
Equivalently, $J_1^{q_1}\leq  J_2^{q_2}$ if and only if for every (some) positive integer $m$
with $mq_1$, $mq_2\in\ZZ$, the ideal $J_1^{mq_1}$ is contained in the integral closure of 
$J_2^{mq_2}$.  It is clear that this is an order relation,
and if $\fra\leq\frb$ and $\frb\leq\fra$, then $\fra=\frb$ in $\Id(X)$. Note also that this order
relation is compatible, in an obvious sense, with sums, products, and rational powers.

\bigskip

We now give the definition of log canonical  pairs in a setting that is general enough for our purpose. Let $X$ be a connected, nonsingular scheme of characteristic zero. We always assume that $X$ is excellent, that is, $X$ has a finite cover by spectra of excellent rings
(see \cite{Matsumura}, Ch. 13 for the definition and basic properties of excellent rings).
The key point of this assumption is that it allows us to use log resolutions, as in the case of schemes of finite type over a field.
In particular, we are in this setting if $X={\rm Spec}(R)$, with
$R=K\[x_1,\ldots,x_n\]$ ($K$ being a field of characteristic zero), or more generally, if $X$
is nonsingular, connected, and of finite type over ${\rm Spec}(R)$.

Let $\fra=J^q$ be a nonzero $\QQ$-ideal on $X$. A log resolution of $(X,\fra)$ is a proper birational 
morphism $\phi\colon X'\to X$, with $X'$ nonsingular, $\fra\cdot\cO_{X'}=\cO(-F)$, with $F$ a divisor, 
and such that the union of $F$ and of the exceptional locus of $\phi$ has simple normal crossings. 
It follows from \cite{Temkin} that there are such log resolutions (and furthermore,
any two log resolutions are dominated by a third one). 
Similarly, it is shown in \emph{loc. cit.}
 that given any $X'\to X$ birational, there is $X''\to X'$
proper, birational, and with $X''$ nonsingular. In particular, given any divisor $E$ over $X$, $E$ appears
as a prime divisor on a nonsingular $X'$, with $X'\to X$ proper and birational. We define
$\ord_E(K_{-/X})$ to be the coefficient of $E$ in the discrepancy divisor $K_{X'/X}$.
It is clear that this depends only on $E$, and not on the model $X'$. 

\begin{definition}
Let $X$ be as above, and $\fra$, $\frb$ nonzero $\QQ$-ideals on $X$.
The pair $(X,\fra\frb^{-1})$ is log canonical if for every divisor $E$ over $X$ we have 
$\ord_E(K_{-/X})+1\geq \ord_E(\fra\frb^{-1})$. 
 \end{definition}
 
 As in the usual setting in \cite{kollar}, one can show that if $X'\to X$ is a log resolution 
 of $(X,\fra\cdot\frb)$, then it is enough to check the defining condition for the prime divisors $E$ on $X'$ (see also 
 \cite{dFM}).  Given a proper, nonzero ideal $J$ in $R$, the \emph{log canonical threshold}
 $\lct(J)$ is defined by
 $$\lct(J):=\sup\{c>0\mid (X,J^c)\, \text{is log canonical}\}.$$ 
 It is easy to get a formula
 for $\lct(J)$ in terms of a log resolution. In particular, this formula shows that $\lct(J)$ is a rational number,
 and that the supremum in the above definition is a maximum (for details see \cite{kollar} and \cite{dFM}).
 Moreover, we have $\lct(J)\leq \dim(X)$, and $\lct(J)\leq 1$ if $J$ is a principal ideal.
We make the convention that $\lct(J)=\infty$ if $J=\cO_X$, and $\lct(J)=0$ if $J=(0)$.

 \bigskip
 
 From now on we will work in the following setting.
  Let $K$ be an algebraically closed 
 field of characteristic zero. 
 Given a positive integer $n$, let $R=K\[x_1,\ldots,x_n\]$, and $X={\rm Spec}(R)$. 
 We sometimes denote by ${\mathbf 0}$ the closed point of $X$.
 The following proposition allows us to check log canonicity only using divisors with center at ${\mathbf 0}$. The argument is well-known (see, for example, Lemma~2.6 in \cite{dFM}).
 
 \begin{proposition}\label{closed_point}
 Given nonzero $\QQ$-ideals $\fra$ and $\frb$ on $X$, the pair $(X,\fra\frb^{-1})$ is log canonical 
 if and only if $\ord_E(K_{-/X})+1\geq \ord_E(\fra\frb^{-1})$ for every divisor $E$ over $X$ 
 with center at ${\mathbf 0}$.
 \end{proposition}
 
 \begin{proof}
 We need to show that the inequality for the divisors with center at ${\mathbf 0}$ implies the inequality for every divisor. Let $\phi\colon X'\to X$ be a log resolution of $(X,\fra\cdot\frb\cdot\frm)$,
 where $\frm$ is the maximal ideal of $R$. 
 Given  a divisor $E$ on $X'$, its image contains the closed point, hence $E$ intersects some component $E_0$ of the inverse image of the closed point (note that by assumption this inverse image is a divisor having simple normal crossings with $E$).

We construct a sequence of divisors $E_m$ over $X$ as follows. Let $E_1$ be the exceptional divisor on the blow-up $X_1$ of $X$ along some connected component of $E\cap E_0$. We define 
$E_m$ and $X_m$ recursively: $X_m$ is the blow-up of some connected component of
the intersection of $E_{m-1}$ with the proper transform of $E$, and $E_m\subset X_m$
is the exceptional divisor of this blow-up. Note that by construction,  each $E_m$ 
has center at ${\mathbf 0}$.
It is standard to check that
$$\ord_{E_m}(\fra\frb^{-1})=m\cdot \ord_{E}(\fra\frb^{-1})+\ord_{E_0}(\fra\frb^{-1}),\,\,\text{and}$$
$$1+\ord_{E_m}(K_{-/X})=m(1+\ord_{E}(K_{-/X}))+(1+\ord_{E_0}(K_{-/X})).$$
Applying the hypothesis for every $E_m$ gives
$$m(1+\ord_E(K_{-/X})-\ord(\fra\frb^{-1}))+(1+\ord_{E_0}(K_{-/X})-\ord_{E_0}(\fra\frb^{-1}))\geq 0$$
for every $m\geq 0$. Therefore $1+\ord_E(K_{-/X})\geq\ord(\fra\frb^{-1})$ for every $E$.
 \end{proof}

 \section{Log canonical pairs via cylinders in the space of arcs}

 We now turn to the characterization of log canonical pairs in terms of arc spaces.
 This characterization is the main tool for our inductive description of singularities. 
 In the case when the ambient space is ${\mathbf A}^n$ (or more generally, a nonsingular variety over $K$), the characterization is given in \cite{Mus} and \cite{ELM}. We need to show that the result extends to our present setting, when the ambient variety is $X={\rm Spec}(R)$, with $R=K\[x_1,\ldots,x_n\]$.
 
 We only consider arcs and jets over the closed point of $X$, so the corresponding spaces can be
 constructed also from the affine space.  Our reference for spaces of arcs and jets for
 schemes of finite type over $K$ is \cite{ELM}. 
 As in the previous section, $K$ is an algebraically closed field of characteristic zero. In addition, 
 we assume that $K$ is uncountable, since 
 since some results on spaces of arcs from \emph{loc. cit.}
 require this assumption.
 
 Given an arbitrary scheme $Y$ over $K$, and
 $y\in Y$ with residue field $K$, the space of arcs on $Y$ centered at $y$ is the set 
 of homomorphisms of local $K$-algebras $\cO_{Y,y}\to K\[t\]$ lifting the morphism to the residue field
 $\cO_{Y,y}\to K$. We denote this space by $(Y,y)_{\infty}$. It is clear from the definition that 
 if $\widehat{Y}={\rm Spec}\,\widehat{\cO_{Y,y}}$, and if $\widehat{y}$ is the closed point of
 $\widehat{Y}$, then we have a canonical identification $(\widehat{Y},\widehat{y})_{\infty}
 =(Y,y)_{\infty}$. In particular, we have $(X,{\mathbf 0})_{\infty}
 =(\AAA^n,{\mathbf 0})_{\infty}$, where we denote by ${\mathbf 0}$ also the origin in the affine space.
 
 Similarly, the space of $m$-jets of $(Y,y)$ is the set of homomorphisms of local $K$-algebras
 $\cO_{Y,y}\to K\[t\]/(t^{m+1})$ that lift the morphism to the residue field. We denote this by $(Y,y)_m$. Note that we have canonical truncation maps 
 $$(Y,y)_{\infty}\to (Y,y)_{m+1}\to (Y,y)_m\to (Y,y)_1=T_yY.$$
 It is clear the we again have identifications $(Y,y)_m=(\widehat{Y}, \widehat{y})_m$. 
 
 From now on we will be interested in the space of arcs $(X,{\mathbf 0})_{\infty}$ that
 we identify as above with $(\AAA^n,{\mathbf 0})_{\infty}$. This space can be further identified
 to $(tK\[t\])^n$, such that an $n$-uple $(u_1,\ldots,u_n)$ corresponds to the morphism 
 $K\[x_1,\ldots,x_n\]\to K\[t\]$ that takes each $x_i$ to $u_i$. We similarly get an identification
 $(X,{\mathbf 0})_m=(tK\[t\]/t^{m+1}K\[t\])^n$, such that the canonical truncation maps
 are given by component-wise
 truncation. We denote by $\psi_m$ the canonical map $(X,{\mathbf 0})_{\infty}\to
 (X,{\mathbf 0})_m$. Note that via the above identifications each $(X,{\mathbf 0})_m$ has an algebraic structure (it is an affine space of dimension $mn$), the truncation maps are morphisms, and 
 $(X,{\mathbf 0})_{\infty}=\projlim_m(X,{\mathbf 0})_m$.
 
 Recall that a cylinder in $(tK\[t\])^n$ is a subset of the form $C=\psi_m^{-1}(S)$ for some $m\geq 1$, and some constructible subset $S$ of $(tK\[t\]/t^{m+1}K\[t\])^n$. We say that $C$ is open, closed,  locally closed, or irreducible  if $S$ is this way (this is independent of the representation of $C$ that we choose). If $S$ is locally closed, then the decomposition of $S$ as a union of irreducible components
 induces a corresponding decomposition for $C$.
  The codimension of a cylinder $C=\psi_m^{-1}(S)$ is defined as 
 $\codim(C):=(m+1)n-\dim(S)$ (despite the fact that we work only in $(\AAA^n,{\mathbf 0})_{\infty}$,
 we follow the convention from \cite{ELM}, computing the codimension in the whole arc space
 $\AAA^n_{\infty}$). 
 
 As in the previous section,
 given $f\in R$ and $\gamma\in (X,{\mathbf 0})_{\infty}$
 that we identify with a morphism $R\to K\[t\]$, the \emph{order of vanishing} 
 $\ord_{\gamma}(f)$ of $f$ along $\gamma$
 is $\ord(\gamma(f))\in\ZZ_{\geq 0}\cup\{\infty\}$. If 
$C\subseteq (X,{\mathbf 0})_{\infty}$ is an irreducible locally closed cylinder, then we put for a nonzero $f\in R$
 $$\ord_C(f):=\min\{\ord_{\gamma}(f)\mid\gamma\in C\}.$$
 It is easy to see that this is finite: if $C=\psi_m^{-1}(S)$, and if $\ord_C(f)=\infty$, then 
 for every 
 $\gamma=(u_1,\ldots,u_n)\in C$, we have $f(u_1+t^{m+1}v_1,\ldots,u_n+t^{m+1}v_n)=0$
 for every $v_1,\ldots,v_n\in K\[t\]$. We deduce using the Taylor expansion that $f=0$.  
 Since we have $\ord_C(f)=\ord_{\gamma}(f)$ for $\gamma$ in a suitable subcylinder of $C$
 that is open in $C$,
 it follows that $\ord_C(fg)=\ord_C(f)+\ord_C(g)$ for every nonzero $f,g\in R$, and therefore
 $\ord_C$ extends (uniquely) to a valuation of the fraction field of $R$, with center at ${\mathbf 0}$
 (this means that $\ord_C(f)\geq 1$ if and only if $f$ lies in the maximal ideal of $R$).
 
 In fact, we have defined more generally $\ord_{\gamma}(\fra)$ for every $\gamma\in
 \left(K\[t\]\right)^n$, and every $\QQ$-ideal $\fra$ in $R$ or in $K[x_1,\ldots,x_n]$. 
 We can define as above $\ord_C(\fra)$ for every
 irreducible locally closed cylinder $C\subseteq (tK\[t\])^n$. If $\fra$ and $\frb$
 are $\QQ$-ideals in $R$,
 then $\fra\leq\frb$ if and only if $\ord_{\gamma}(\fra)\geq\ord_{\gamma}(\frb)$ for every
 $\gamma\in (tK\[t\])^n$. A similar assertion holds for $\QQ$-ideals in $K[x_1,\ldots,x_n]$
 if we evaluate along all $\gamma\in \left(K\[t\]\right)^n$.
 
 A fundamental example of cylinders is given by the contact loci of a given ideal. If $\fra$ is a
 nonzero $\QQ$-ideal in $R$, then ${\rm Cont}^{m}(\fra):=\{\gamma\mid\ord_{\gamma}(\fra)=m\}$. We similarly define ${\rm Cont}^{\geq m}(\fra)$. It is clear that ${\rm Cont}^{\geq m}(\fra)$ is a 
 closed cylinder, while ${\rm Cont}^m(\fra)$ is a locally closed cylinder.

 We now explain how to relate divisors and cylinders via valuation theory, following \cite{ELM}. 
 The key point is that divisors over $X$ centered at ${\mathbf 0}$ can be identified with divisors over
 $\AAA^n$ with center at ${\mathbf 0}$. More precisely: if $E$ is a divisor centered at ${\mathbf 0}$ over $\AAA^n$, 
 then $E$ is a prime divisor on some $Y$, having a proper, birational morphism to $\AAA^n$
 that maps $E$ to ${\mathbf 0}$. Base-changing via $X\to\AAA^n$ gives a prime divisor
 $\widehat{E}$ on $X':=X\times_{\AAA^n}Y$ lying over ${\mathbf 0}\in X$. This is the unique divisor 
 $\widehat{E}$ over $X$ such that $\ord_{\widehat{E}}$ agrees with $\ord_E$ on $\cO_{\AAA^n,{\mathbf 0}}$.  Furthermore, every divisor over $X$ with center at ${\mathbf 0}$ can be written as
 $\widehat{E}$ for a unique such $E$.
 
We note that more generally, there is a bijection between valuations of the function field of
$\AAA^n$ with center at ${\mathbf 0}$, and valuations of the function field of $X$ with center at ${\mathbf 0}$
given by restriction from $R$ to $\cO_{\AAA^n,{\mathbf 0}}$
(this is a bijection due to the fact that such valuations are automatically continuous for the 
$(x_1,\ldots,x_n)$-adic topology on the respective rings). Of course, via this bijection 
$\ord_{\widehat{E}}$ corresponds to $\ord_E$, and if $C$ is an irreducible closed cylinder in 
$(X,{\mathbf 0})_{\infty}$, then $\ord_C$ corresponds to the valuation associated to $C$ when considered as a cylinder in $(\AAA^n,{\mathbf 0})_{\infty}$. 

 Given a divisor $E$ over $\AAA^n$ with center ${\mathbf 0}$, for every positive integer $m$
 one can associate an irreducible closed cylinder ${\mathcal C}_m(E)\subseteq (tK\[t\])^n$, see 
 \cite{ELM}. If $E$ is a prime divisor on some variety $Y$ over $\AAA^n$, then 
 ${\mathcal C}_m(E)$ is the closure of the image of ${\rm Cont}^m(E)\subseteq Y_{\infty}$
 in $(\AAA^n,{\mathbf 0})_{\infty}$. An easy fact that follows from definition is that 
 $\ord_{{\mathcal C}_m(E)}=m\cdot\ord_E$. It is shown in Theorem~2.1 in \emph{loc. cit.}  that
 $\codim({\mathcal C}_m(E))=
 m\cdot (\ord_E(K_{-/\AAA^n})+1)$. 
 Furthermore, given any irreducible closed cylinder $C$ in $(\AAA^n,{\mathbf 0})_{\infty}$, there is 
 a divisor $E$ over $\AAA^n$ with center at the origin, and a positive integer $m$ such that
 $C\subseteq {\mathcal C}_m(E)$ and
 $\ord_C=\ord_{{\mathcal C}_m(E)}$ (see Theorem~2.7 and its proof in \emph{ibid}). By the previous discussion, we can translate these results as a correspondence between divisors over $X$ with center at ${\mathbf 0}$, and irreducible
 closed cylinders in $(X,{\mathbf 0})_{\infty}$. 
 
  We now can prove the interpretation of log canonicity of a pair $(X,\fra\cdot\frb^{-1})$
 in terms of the space of arcs. The proof is a small variation on the argument in \cite{ELM}.

 \begin{theorem}\label{description}
Suppose that $K$ is uncountable, and let $\fra$ and $\frb$ be nonzero $\QQ$-ideals in $R$.
\begin{enumerate}
\item[i)] The pair $(X,\fra\cdot\frb^{-1})$ is log canonical if and only if for every irreducible closed cylinder 
$C$ in $(X,{\mathbf 0})_{\infty}$ we have
$$\ord_C(\frb)\geq\ord_C(\fra)-\codim(C).$$
\item[ii)] Given $P_1,\ldots,P_s\in R$, and a positive integer $m$, it is enough to check the condition in ${\rm i)}$ only for those $C$ with the property that $m\vert \ord_{C}(P_i)$ for every $i$.
\end{enumerate}
\end{theorem}

\begin{proof}
Let $m$ be a positive integer, and $E$ a divisor over $X$ with center ${\mathbf 0}$. Since
$\ord_{{\mathcal C}_m(E)}=m\cdot \ord_E$, it follows that
$$\ord_{{\mathcal C}_m(E)}(\fra)=m\cdot \ord_E(\fra),\,\,\text{and}\,\,
\ord_{{\mathcal C}_m(E)}(\frb)=m\cdot\ord_E(\frb).$$
On the other hand, $\codim({\mathcal C}_m(E))=m(\ord_E(K_{-/X})+1)$. 
Therefore the condition in i) for $C={\mathcal C}_m(E)$ (for some $m$) is equivalent to the condition for the log canonicity of $(X,\fra\frb^{-1})$ for the divisor $E$. 
We deduce from Proposition~\ref{closed_point} that the condition in i) for every irreducible closed cylinder $C$ in $(X,{\mathbf 0})_{\infty}$ implies that
$(X,\fra\frb^{-1})$ is log canonical. Since it is enough to use $C={\mathcal C}_m(E)$, with
$m$ divisible enough, in which case $\ord_C(P_i)$ also is divisible enough, we get the assertion in ii). 

The converse in i) follows from the fact that given an arbitrary closed,
irreducible cylinder $C$, there is a divisor $E$ with center ${\mathbf 0}$, and a positive integer $m$
such that $C\subseteq {\mathcal C}_m(E)$, and the two cylinders define the same valuation. Therefore
$$\ord_C(\fra)=m\cdot\ord_E(\fra),\,\,\ord_C(\frb)=m\cdot\ord_E(\frb),\,\,\text{and}$$
$$\codim(C)\geq\codim({\mathcal C}_m(E))=m(1+\ord_E(K_{-/X})).$$
This completes the proof of i). 
\end{proof}

\begin{remark}
The condition in Theorem~\ref{description} that $C$ be closed is not essential: given any irreducible locally closed cylinder $C$, the condition in i) for $C$ is the same as the condition for $\overline{C}$.
It follows from part ii) in the theorem that given $P_1,\ldots,P_r\in R$, in order to show that $(X,\fra\cdot\frb^{-1})$ is log canonical, it is enough to check the condition
in i) only for those irreducible locally closed cylinders $C$ with the property that $m\vert
\ord_{\gamma}(P_i)$ for every $i$, and every $\gamma\in C$.
\end{remark}

\section{Orders of symmetric polynomials in the roots}

In this section $K$ denotes an algebraically closed field of characteristic zero.
We consider  polynomials
$h=y^d+a_1y^{d-1}+\cdots+a_d$, with $a_i\in K\[t\]$ for all $i$.
Let $\alpha_1,\ldots,\alpha_d$ be the roots of $h$. The
$\alpha_i$ are Puiseux formal power series, lying in $K\[t^{1/N}\]$,
where we may take $N=d!$. 
We denote by $s_i$ the $i^{\rm th}$
elementary symmetric polynomial. Therefore
$a_i=(-1)^is_i(\alpha_1,\ldots,\alpha_d)$.

\begin{lemma}\label{lem1}
If $\frb$ is the $\QQ$-ideal $\sum_{i=1}^d(z_i)^{1/i}$ in $\QQ[z_1,\ldots,z_d]$, then
for every $h$ as above we have
$$\min_i\{\ord(\alpha_i)\}=\ord(\frb(a_1,\ldots,a_d)).$$
\end{lemma}

\begin{proof}
We need to show that for every $q\in\frac{1}{N}\ZZ_{\geq 0}$, we have
$\ord(\alpha_i)\geq q$ for every $i$, if and only if $\ord(a_i)\geq iq$ for every $i$.
Since $a_i=(-1)^{i}s_i(\alpha_1,\ldots,\alpha_d)$, it is clear
that if $\ord(\alpha_i)\geq q$ for every $i$, then $\ord(a_i)\geq
iq$ for every $i$. Conversely, if $\ord(a_i)\geq iq$ for every
$i$, then $h':=y^d+\sum_{i=1}^dt^{-iq}a_iy^{d-i}$
has coefficients in $K\[t^{1/N}\]$. For every $i$ we have $h(\alpha_i)=0$, and since 
$K\[t^{1/N}\]$ is normal, we deduce that
$t^{-q}\alpha_i\in K\[t^{1/N}\]$.
\end{proof}

\begin{remark}
Note that for the equality in the lemma we do not need to assume that the two 
quantities are finite.
\end{remark}

We now fix $k\in\{1,\ldots,d\}$. If we consider all products
$\alpha_{j_1}\cdot\ldots\cdot\alpha_{j_k}$, where $j_1,\ldots,j_k$ are
pairwise distinct, then applying $s_{\ell}$ to these products gives
a symmetric polynomial in $\alpha_1,\ldots,\alpha_d$.
Therefore we can write it as $A_k^{(\ell)}(a_1,\ldots,a_d)$, for some $A_k^{(\ell)}
\in\QQ[z_1,\ldots,z_d]$, independent of the ground field $K$. 
Consider the $\QQ$-ideal $\frb_k:=\sum_{\ell=1}^{{d}\choose{k}}\left(A_k^{(\ell)}\right)^{1/\ell}$
in $\QQ[z_1,\ldots,z_d]$.
Lemma~\ref{lem1} implies

\begin{corollary}\label{cor1}
For every $K$ and $h$ as above, and for every $k\leq d$, we have
$$\min\{\ord(\alpha_{j_1}\cdot\ldots\cdot\alpha_{j_k})\}
=\ord(\frb_k(a_1,\ldots,a_d)),$$ where the minimum on the left-hand side is over
all distinct $j_1,\ldots,j_k\in\{1,\ldots,d\}$.
\end{corollary}

\begin{remark}
The above procedure is not always the most economical. For example, when $k=d$
we see right away that we can take $\frb_d=(z_d)$. 
\end{remark}

We are also interested in the following situation:  consider $\{\alpha_{\ell}-w\mid 1\leq
\ell\leq d\}$, for some $w\in K\[t^{1/N}\]$. Note that these are precisely the roots of
$$h(y+w)= y^d+\sum_{j=1}^d\widetilde{a}_jy^{d-j},$$
where $\widetilde{a}_j= \frac{1}{(d-j)!}\cdot\frac{\partial^{d-j}
h}{\partial y^{d-j}}(w)$. Given $k$, we can write
$$A_k^{(\ell)}(\widetilde{a}_1,\ldots,\widetilde{a}_d)=\widetilde{A}_k^{(\ell)}(a_1,\ldots,a_d,w),$$
for a unique polynomial $\widetilde{A}_k^{(\ell)}\in\QQ[z_1,\ldots,z_{d+1}]$.
If we denote by $\widetilde{\frb}_k$ the $\QQ$-ideal
$\sum_{\ell=1}^{{d}\choose{k}}\left(\widetilde{A}_k^{(\ell)}\right)^{1/\ell}$, then
Corollary~\ref{cor1} implies

\begin{corollary}\label{cor2}
For every $K$, $h$, and $w$ as above, and for every $k\leq d$, we have
$$\min\{\ord(\alpha_{j_1}-w)\cdots(\alpha_{j_k}-w)\}
=\ord(\widetilde{\frb}_k(a_1,\ldots,a_d,w)),$$ where the minimum on the
left-hand side is over all distinct $j_1,\ldots,j_k\in\{1,\ldots,d\}$.
\end{corollary}

We now use the above ideas to give a criterion for a polynomial to
have all its roots in $K\[t\]$.

\begin{corollary}\label{cor3}
Given $d$, there are polynomials $P_1,\ldots,P_r
\in\QQ[z_1,\ldots,z_d]$ and a positive integer $m$ such
that for every $K$ and every polynomial $h=y^d+\sum_{i=1}^da_iy^{d-i}\in K\[t\][y]$, if $m$
divides $\ord(P_i(a_1,\ldots,a_d))$ for every $i$, then all roots
$\alpha_1,\ldots,\alpha_d$ of $h$ are in $K\[t\]$.
\end{corollary}

\begin{proof}
Let $\alpha_1,\ldots,\alpha_d$ be the roots of $h$. Suppose that
$\alpha_i$ is a root that lies in $K\[t^{1/N}\]$, but not
in any $K\[t^{1/m}\]$, with $m$ dividing $N$. If we write
$\alpha_i=h(t^{1/N})$, and if $\eta$ is a primitive root of
order $N$ of $1$, then $h(\eta t^{1/N})$ is equal to $\alpha_j$, for some
$j\neq i$. In this case $\ord(\alpha_i-\alpha_j)$ is finite, but it
is not an integer. Therefore all $\alpha_i$ are in $K\[t\]$ if and only if all
$\ord(\alpha_i-\alpha_j)$ are in $\ZZ\cup\{\infty\}$.

Let $H$ be the polynomial whose roots are the $(\alpha_i-\alpha_j)$,
with $i\neq j$. It is clear that we may write
$$H(y)=y^D+\sum_{i=1}^Db_i(a_1,\ldots,a_d)y^{D-i},$$
for some polynomials $b_i\in\QQ[z_1,\ldots,z_d]$, with
$1\leq i \leq D={{d}\choose{2}}$. Let us order the roots
$\beta_1,\ldots,\beta_D$ of $H$ such that
$$\ord(\beta_1)\leq\ldots\leq\ord(\beta_D)\leq\infty.$$
It follows from Corollary~\ref{cor1} that we can find $\QQ$-ideals
$\frc_k$ in $\QQ[w_1,\ldots,w_D]$ such that
$$\ord(\beta_1)+\ldots+\ord(\beta_k)=\ord(\frc_k(b_1(a),\ldots,b_D(a))).$$ 
For some positive integer $m$, we can write
$\frc_k=(Q_{k,1}(w),\ldots,Q_{k,r_k}(w))^{1/m}$ 
for every $k$.
We see that $m$, and the set consisting of all the polynomials
$Q_{i,j}(b_1,\ldots,b_D)$ satisfy the 
condition in the corollary.
\end{proof}

\bigskip

We give below a different way of computing the minimum in
Corollary~\ref{cor1}. This can be useful in practice, having the
advantage of being more explicit. It is essentially the expression
that comes out of the Newton method (see \cite{koblitz}, \S IV.3).

Given $k\leq d$, and not necessarily distinct 
$i_1,\ldots,i_k\in\{1,\ldots,d\}$ with $i_m\geq k-m+1$ for every $m\in\{1,\ldots,k\}$, we 
consider the $\QQ$-ideal 
$M(i_1,\ldots,i_k):=\prod_{m=1}^k \left(z_{i_m}\right)^{j_m}$ in 
$\QQ[z_1,\ldots,z_d]$, where
$j_m={\frac{1}{i_m-k+m}\cdot
\prod_{\ell=1}^{m-1}\frac{i_{\ell}-k+\ell-1}{i_{\ell}-k+\ell}}$
(note that this $\QQ$-ideal is a rational power of a principal monomial ideal).
We put $\overline{\frb}_k:=\sum M(i_1,\ldots,i_k)$, the sum being over all
$(i_1,\ldots,i_k)$ as above. Note that $\overline{\frb}_1$ is equal to the
$\QQ$-ideal $\frb$ defined in Lemma~\ref{lem1}.

\begin{proposition}\label{prop1}
Let $h=y^d+\sum_{i=1}^da_iy^{d-i}\in K\[t\][y]$, having the roots
$\alpha_1,\ldots,\alpha_d$. For every $k$ with $1\leq k\leq d$, we
have
$$\min\{\ord(\alpha_{j_1}\cdot\ldots\cdot\alpha_{j_k})\}=
\ord(\overline{\frb}_k(a_1,\ldots,a_d)),$$
where the minimum on the
left-hand side is over all pairwise distinct $j_1,\ldots,j_k\in\{1,\ldots,d\}$.
\end{proposition}

\begin{proof}
Let $q_i=\ord(\alpha_i)$. We may assume that the roots of $h$ are
labeled such that $q_1\leq\ldots\leq q_d$.
 We prove the assertion by induction on $k$, the case $k=1$ being
covered by Lemma~\ref{lem1}. In order to prove the induction step,
it is enough to show the following equality
\begin{equation}\label{eq1}
\sum_{j=1}^kq_j=\min\left\{\frac{\ord(a_i)}{i-k+1}+\frac{i-k}{i-k+1}\cdot\sum_{j=1}^{k-1}q_j\vert
k\leq i\leq d \right\}.
\end{equation}
For (\ref{eq1}), note that if $k\leq i\leq d$, then
\begin{equation}\label{eq2}
q_1+\ldots+q_{k-1}+(i-k+1)q_k\leq \sum_{j=1}^iq_j\leq\ord(a_i).
\end{equation}
Moreover, if $i$ is chosen such that $q_k=\ldots=q_{i}<q_{i+1}$ (if $q_k=q_d$, then we take
$i=d$), then both inequalities in (\ref{eq2}) become
equalities. These two facts imply the formula (\ref{eq1}).
\end{proof}

In the case $k=d-1$, one can use directly the formula in
Lemma~\ref{lem1} to obtain a more compact expression. This in turn
can be used to describe the maximum of the order of the roots of
$f$.

\begin{corollary}\label{cor4}
If $\frc$ is the $\QQ$-ideal 
$\sum_{i=1}^d(z_{d-i})^{1/i}(z_d)^{(i-1)/i}$ in $\QQ[z_1,\ldots,z_d]$, then
for every polynomial $h=y^d+\sum_{i=1}^da_iy^{d-i}\in K\[t\][y]$ as above, with
roots $\alpha_1,\ldots,\alpha_d$, we have
\begin{enumerate}
\item[i)] $\min_i\left\{\sum_{j\neq i}\ord(\alpha_j)\right\}=\ord(\frc(a_1,\ldots,a_d))$.
\item[ii)] $\max_i\{\ord(\alpha_i)\}=\ord(a_d)-\ord(\frc(a_1,\dots,a_d))$.
\end{enumerate}
\end{corollary}

\begin{proof}
Let $\beta_{\ell}:=\prod_{j\neq
\ell}\alpha_j=(-1)^d\frac{a_d}{\alpha_{\ell}}$. Note that we have
$s_i(\beta_1,\ldots,\beta_d)=\pm (a_d)^{i-1}a_{d-i}$
for $1\leq i\leq d$, where we put $a_0=1$.
The formula in i) is then a consequence of
Lemma~\ref{lem1}. Note that the expression in i) is equal to
$\ord(\alpha_1\cdot\ldots\cdot\alpha_d)-\max_i\{\ord(\alpha_i)\}$,
hence we get the formula in ii).

\end{proof}

\bigskip

In the next section we will make use of 
the following proposition that shows how to take the maximum 
of expressions as in Corollary~\ref{cor2}, when we let $w$ vary over all $\alpha_i$.

\begin{proposition}\label{lem1_2}
Given a $\QQ$-ideal $\frb$  in $\QQ[z_1,\ldots, z_{d+1}]$, there are $\QQ$-ideals 
$\frb_+$ and $\frb_-$ in $\QQ[z_1,\ldots,z_d]$ such that for every $K$ and every $h=y^d+\sum_{i=1}^d
a_iy^{d-i}\in K\[t\][y]$, with roots $\alpha_1,\ldots,\alpha_d$, we have
$$\max\{\ord(\frb(a_1,\ldots,a_d,\alpha_i)\mid 1\leq i\leq d\}=\ord(\frb_+(a_1,\ldots,a_d))
-\ord(\frb_-(a_1,\ldots,a_d))$$
${\rm (}$we follow the convention that $\infty-\infty=\infty$${\rm )}$.
\end{proposition} 

\begin{proof}
We use the idea in the proof of Corollary~\ref{cor4}. For every
$k\in\{1,\ldots,d\}$ we put
$$\lambda_k(a_1,\ldots,a_d)=\min_{j_1,\ldots,j_k}\{\ord\prod_{\ell=1}^k\frb(a_1,\ldots,a_d,\alpha_{j_{\ell}})\},$$
the minimum being over all distinct $j_1,\ldots,j_k\in\{1,\ldots,d\}$.
Since
$$\max_i\{\ord(\frb(a_1,\ldots,a_d,\alpha_i)\mid 1\leq i\leq d\}=\lambda_d(a_1,\ldots,a_d)-
\lambda_{d-1}(a_1,\ldots,a_d),$$
we see that in order to prove the proposition it is enough to show that for every $k$, there is a 
$\QQ$-ideal $\frc_k$ in $\QQ[z_1,\ldots,z_d]$ such that
$\lambda_k(a_1,\ldots,a_d)=\ord(\frc_k(a_1,\ldots,a_d))$ for every 
$a_1,\ldots,a_d\in K\[ t\]$. 

Let us put $a=(a_1,\ldots,a_d)$. If $\frb=(P_1,\ldots,P_r)^q$, we see that
$$\lambda_k(a_1,\ldots,a_d)=q^k\cdot\min\{\ord(P_{i_1}(a,\alpha_{j_1})\cdot\ldots\cdot
P_{i_k}(a,\alpha_{j_k}))\},$$
where the minimum is over the distinct $j_1,\ldots,j_k\in\{1,\ldots,d\}$, and over the
not necessarily distinct $i_1,\ldots,i_k$ in the same set. Note that every symmetric function of
 all the $P_{i_1}(a,\alpha_{j_1})\cdot\ldots\cdot P_{i_k}(a,\alpha_{j_k})$
can be written as a polynomial with rational coefficients in $a_1,\ldots,a_d$ and the symmetric functions in $\alpha_1,\ldots,
\alpha_d$, hence just as a polynomial in the $a_i$'s. The existence of $\frc_k$ is now a consequence of Lemma~\ref{lem1}.
\end{proof}

\begin{remark}\label{containment}
It follows from the above proof that we can find $\frb_+$ and $\frb_-$ as in the statement of
Proposition~\ref{lem1_2} such that for every $K$ we have 
$$\frb_+ K[z_1,\ldots,z_d]\leq \left(\frb_- K[z_1,\ldots,z_d]\right)^{\frac{d}{d-1}}.$$ 
Indeed, for this it is enough
to show that for every $a_1,\ldots,a_d\in K\[t\]$, we have
$$\ord(\frb_+(a_1,\ldots,a_d))\geq\frac{d}{d-1}\cdot\ord(\frb_-(a_1,\ldots,a_d).$$
In order to see this, note that if $\alpha_1,\ldots,\alpha_d$ are the roots of
$y^d+\sum_{i=1}^da_iy^{d-i}$, and if we label the $\alpha_i$ such that
$$\ord(\frb(a_1,\ldots,a_d,\alpha_1))\leq\ldots\leq\ord(\frb(a_1,\ldots,a_d,\alpha_d)),$$
then $\lambda_{d-1}(a_1,\ldots,a_d)=\sum_{i=1}^{d-1}\ord(\frb(a_1,\ldots,a_d,\alpha_i))
\leq (d-1)\ord(\frb(a_1,\ldots,a_d,\alpha_d))$, which implies 
$$\lambda_{d}(a_1,\ldots,a_d)=\lambda_{d-1}(a_1,\dots,a_d)+
\ord(\frb(a_1,\ldots,a_d,\alpha_d))\geq\frac{d}{d-1}\cdot \lambda_{d-1}(a_1,\dots,a_d).$$
This gives our statement.
\end{remark}

We now deduce the following corollary, that will be applied in the next section to get the inductive description of singularities of pairs.

\begin{corollary}\label{cor_section2}
For every positive integers $d$ and $k$, with $k\leq d$, and every nonnegative rational numbers $c_1$, $c_2$, there are $\QQ$-ideals 
$\frp_+=\frp_+(d,k,c_1,c_2)$ and $\frp_-=\frp_-(d,k,c_1,c_2)$ in $\QQ[z_1,\ldots,z_d]$ such that 
for every $K$ we have
\item[i)] $\frp_+ K[x_1,\ldots,x_d]\leq \left(\frp_- K[x_1,\ldots,x_d]\right)^{\frac{d}{d-1}}$.
\item[ii)] For all
 $(a_1,\ldots,a_d)\in (K\[t\])^d$, if $\alpha_1,\ldots,\alpha_d$ are the roots of
$y^d+\sum_{i=1}^da_iy^{d-i}$, then
$$\max_{i=1}^d\left\{c_1\cdot \min_{j_1,\ldots,j_{k-1}}\ord\left(\prod_{\ell=1}^{k-1}(\alpha_{j_{\ell}}-\alpha_i)\right)+c_2\cdot \min_{j_1,\ldots,j_k}\ord\left(\prod_{\ell=1}^k(\alpha_{j_{\ell}}-\alpha_i)\right)
\right\}$$
$$=\ord(\frp_+(a_1,\ldots,a_d))-\ord(\frp_-(a_1,\ldots,a_d)),$$
where each minimum on the left-hand side is over distinct elements of $\{1,\ldots,d\}$.
\end{corollary}

\begin{proof}
Using the notation in Corollary~\ref{cor2}, we 
define the following $\QQ$-ideal in $\QQ[z_1,\ldots,z_d]$ 
$$\frp:=(\widetilde{\frb}_{k-1})^{c_1}\cdot(\widetilde{\frb}_k)^{c_2}.$$
It follows that $\max_{i}\{\ord(\frp(a_1,\ldots,a_d))\}$
is equal to the left-hand side of the formula in the statement. 
Applying Proposition~\ref{lem1_2} and Remark~\ref{containment} for $\frp$, we get the $\QQ$-ideals 
$\frp_+$ and $\frp_-$ with the required property.
\end{proof}

\begin{remark}\label{vanishing}
If in the corollary we write $\frp_+=J^q$, then $J\subseteq (x_1,\ldots,x_d)$, unless $c_2=0$, and either
$k=1$, or $c_1=0$.
Indeed, if $a_i\in tK\[t\]$ for all $i$, then $\ord(\alpha_j)>0$ for every $j$, by Lemma~\ref{lem1}. 
If $J\not\subseteq 
(x_1,\ldots,x_d)$, then the expression on the right-hand side in ii) is non-positive,
and the one on the left-hand side is positive, unless $c_2=0$, and 
either $c_1=0$ or $k=1$.
\end{remark}

\begin{example}\label{example0}
With the notation in Corollary~\ref{cor_section2}, suppose that $k=1$. We claim that if $a_1=0$, then
for every $i$ we have $\min_j\{\ord(\alpha_j-\alpha_i)\}=\min_j\{\ord(\alpha_j)\}$. Indeed, the inequality
``$\geq$" is clear, and the opposite inequality follows from
$\sum_{j=1}^d(\alpha_j-\alpha_i)=-d\alpha_i$. We deduce from Lemma~\ref{lem1} that if $a_1=0$, then
$$\max_i\left\{c\cdot \min_j\{\ord(\alpha_j-\alpha_i)\}\right\}=\ord(\frb^c(a_2,\ldots,a_d)),$$
where $\frb=\sum_{j=2}^d(z_j)^{1/j}$, and $c$ is a positive rational number.
\end{example}

\begin{example}\label{example1}
Let $h=y^3+ay+b$, with $a,b\in K\[t\]$, and let
$\alpha_1,\alpha_2$, and $\alpha_3$ be the roots of $h$. After
relabeling the roots, we may assume that
$$\ord(\alpha_1-\alpha_2)=\ord(\alpha_1-\alpha_3)\leq\ord(\alpha_2-\alpha_3).$$
If we fix $i$, then the $(\alpha_j-\alpha_i)$, with $1\leq j\leq 3$, are the roots of
$$h(y+\alpha_i)=y^3+3\alpha_iy^2+(a+3\alpha_i^2)y.$$
Therefore $\ord\prod_{j\neq
i}(\alpha_j-\alpha_i)=\ord(a+3\alpha_i^2)$.

It is a straightforward computation to show that the
$(a+3\alpha_i^2)$, with $1\leq i\leq 3$, are the roots of
$$y^3+3ay^2-\Delta_h=0,$$
where $\Delta_h$ is the discriminant of $h$, that is,
$\Delta_h=4a^3+27b^2$. 
Since
$\Delta_h=\pm(\alpha_1-\alpha_2)^2(\alpha_2-\alpha_3)^2(\alpha_3-\alpha_1)^2$,
we deduce that
$$\ord(\alpha_2-\alpha_3)=\max_{i\neq j}\left\{\ord(\alpha_i-\alpha_j)\right\}=
\ord\prod_{i\neq j}(\alpha_i-\alpha_j)-\min_{i}\{\ord(a+3\alpha_i^2)\}$$
$$=\frac{1}{2}\ord(\Delta_h)-\min\left\{\ord(a),\frac{1}{3}\ord(\Delta_h)\right\}
=\ord((\Delta_h)^{1/2})-\ord((a^3,b^2)^{1/3}).$$

Since
$2\ord(\alpha_1-\alpha_2)+\ord(\alpha_2-\alpha_3)=\frac{1}{2}\ord(\Delta_h)$,
we see that
$$\ord(\alpha_1-\alpha_2)=\ord(\alpha_1-\alpha_3)=\ord((a^3,b^2)^{1/6}).$$

Let us consider for $h$ the situation in Corollary~\ref{cor_section2}, with $k=2$.
In this case 
$$\max_{i=1,2,3}\left\{c_1\cdot \min_{j\neq i}\ord(\alpha_j-\alpha_i)+
c_2\cdot\ord\left(\prod_{j\neq i}(\alpha_j-\alpha_i)\right)\right\}=
c_1\ord(\alpha_1-\alpha_2)$$
$$+c_2\ord\left((\alpha_1-\alpha_2)(\alpha_2-\alpha_3)\right)=
\ord(\frp_+(a,b))-\ord(\frp_-(a,b)),$$
where the $\QQ$-ideals $\frp_+$ and $\frp_-$ in $\QQ[u,v]$ are given by
$$\frp_+= \left\{
\begin{array}{cl}
(4u^3+27v^2)^{\frac{c_2}{2}}, & {\rm if}\,c_2\geq c_1; \\[2mm]
(4u^3+27v^2)^{\frac{c_2}{2}}\cdot (u^3,v^2)^{\frac{c_1-c_2}{6}}, & {\rm if}\,c_2\leq c_1.
\end{array}\right.
$$

$$\frp_- = \left\{
\begin{array}{cl}
(u^3,v^2)^{\frac{c_2-c_1}{6}}, & {\rm if}\,c_2\geq c_1; \\[2mm]
\QQ[u,v], & {\rm if}\,c_2\leq c_1.
\end{array}\right.
$$
Note that in this case we preferred to follow a slightly simplified procedure
for obtaining $\frp_+$ and $\frp_-$, from the one described in
the proof of Corollary~\ref{cor_section2}.

When $k=1$, the situation is simpler: we have
$$\max_{i=1,2,3}\left\{c_2\cdot \min_{j\neq i}\ord(\alpha_j-\alpha_i)\right\}=
c_2\ord(\alpha_1-\alpha_2),$$
and we may take $\frp_+=(u^3,v^2)^{c_2/6}$, and $\frp_{-}=\QQ[u,v]$.
\end{example}

For future reference, we state one more consequence of Corollary~\ref{cor2}.
Given any $w\in K\[t^{1/N}\]$, and any $h=y^d+\sum_{i=1}^da_iy^{d-i}$, we 
put $b_i=\ord(\alpha_i-w)$, where the roots $\alpha_i$ of $h$ are labeled so that
$$\ord(\alpha_1-w)\leq\ldots\leq\ord(\alpha_d-w)=\infty.$$

\begin{corollary}\label{cor100}
For every $k\leq d$, there are $\QQ$-ideals $\frp^{(k)}_+$ and $\frp^{(k)}_-$
 in $\QQ[z_1,\ldots,z_{d+1}]$, such that 
for every $K$,  every $a_1,\ldots,a_d\in K\[t\]$, 
and every $w\in K\[t^{1/N}\]$, 
we have
$$\ord(\alpha_k-w)=\ord(\frp^{(k)}_+(a_1,\ldots,a_d,w))-\ord(\frp^{(k)}_-(a_1,\ldots,a_d,
w)),$$
where $\alpha_1,\ldots,\alpha_d$ are the roots of $y^d+\sum_{i=1}^da_iy^{d-i}$, ordered as described above.
\end{corollary}

\begin{proof}
By the definition of $b_k$ we have
$$b_k=\min_{j_1,\ldots,j_k}\{\ord(\alpha_{j_1}-w)\cdots(\alpha_{j_k}-w)\}-
\min_{j_1,\ldots,j_{k-1}}\{\ord(\alpha_{j_1}-w)\cdots(\alpha_{j_{k-1}}-w)\},$$
where both minima are over distinct elements of $\{1,\ldots,d\}$. 
With the notation in Corollary~\ref{cor2}, we see that we may take
$\frp^{(k)}_+=\widetilde{\frb}_k$, and $\frp^{(k)}_-=\widetilde{\frb}_{k-1}$.
\end{proof}

We end this section with two more propositions that will be needed in the proof of our main results. 
In both statements we consider an irreducible locally closed cylinder $T\subseteq \left(tK\[t\]\right)^d$, 
and $\QQ$-ideals
$\frb$, 
$\frq_1,\ldots,\frq_r$ in $\QQ[z_1,\ldots,z_{d+1}]$. For every 
$a=(a_1,\ldots,a_d)\in T$, we consider the roots of
$y^d+\sum_{\ell=1}^da_{\ell}y^{d-\ell}$, that we denote by $\alpha_1(a),\ldots,\alpha_d(a)$
to emphasize the dependence on $a$. We assume that these roots are in $K\[t\]$
(hence in $tK\[t\]$, by Lemma~\ref{lem1}) for every $a\in T$.

\begin{proposition}\label{subcylinder1}
With the above notation, if there is $a\in T$ with 
$$\max_j\{\ord(\frb(a_1,\ldots,a_d,\alpha_j(a)))\}<\infty,$$
then
 there is a subcylinder $T^{\circ}$ of $T$ that is open in $T$, 
and $N_1,\ldots,N_r\in \ZZ_{\geq 0}\cup\{\infty\}$ such that
for every $a=(a_1,\dots,a_d)\in T^{\circ}$, there
is $i$ with 
$$\ord(\frb(a_1,\ldots,a_d,\alpha_i(a)))=\max_j\{\ord(\frb(a_1,\ldots,a_d,\alpha_j(a)))\},$$
and $\ord(\frq_{\ell}(a_1,\ldots,a_d,\alpha_i(a)))=N_{\ell}$ for $\ell=1,\ldots,r$.
\end{proposition}

\begin{proposition}\label{subcylinder2}
With the above notation, if we have for every $a\in T$
$$\max_j\{\ord(\frb(a_1,\ldots,a_d,\alpha_j(a)))\}=\infty,$$ then 
for every positive integer $L$ there is a subcylinder $T^{\circ}$ of $T$
that is open in $T$, 
and $N_1,\ldots,N_r\in \ZZ_{\geq 0}\cup\{\infty\}$ such that
for every $a=(a_1,\dots,a_d)\in T^{\circ}$, there
is $i$ with $\ord(\frb(a_1,\ldots,a_d,\alpha_i(a)))\geq L$
and $\ord(\frq_{\ell}(a_1,\ldots,a_d,\alpha_i(a)))=N_{\ell}$ for $\ell=1,\ldots,r$.
\end{proposition}

Before proving the propositions, we make some preparations.
Recall the following well-known fact.

\begin{lemma}\label{lem11}
Let $f=y^d+\sum_{i=1}^da_iy^{d-i}$, and $g=y^d+\sum_{i=1}^db_iy^{d-i}$, with all $a_i$, $b_i\in
K\[t\]$, and such that $\ord(a_i-b_i)\geq N$ for aIl $i$. If  $\alpha_1,\ldots,\alpha_d$
are the roots of $f$, and if $\beta$ is a root of $g$, then there is $i$ such that 
$\ord(\beta-\alpha_i)\geq\frac{N}{d}$.
\end{lemma}

\begin{proof}
For every $i$, the order of 
$$\prod_{i=1}^d(\beta-\alpha_i)=f(\beta)=f(\beta)-g(\beta)=\sum_{j=1}^d(a_j-b_j)\beta^{d-j}$$
is $\geq N$, hence $\ord(\beta-\alpha_i)\geq\frac{N}{d}$ for 
some $i$.
\end{proof}

We introduce some notation that will be used in the proofs of Propositions~\ref{subcylinder1}
and \ref{subcylinder2}. Let 
$$\Phi\colon \left(tK\[t\]\right)^d\to \left(tK\[t\]\right)^d,\,\,
\Phi(\alpha_1,\ldots,\alpha_d)=(-s_1(\alpha),\ldots,(-1)^ds_d(\alpha)),$$ where
$s_i$ is the $i^{\rm th}$ elementary symmetric polynomial. If $S=\Phi^{-1}(T)$, then $S$ is a locally closed cylinder,
and our assumption implies that the induced map $S\to T$ is surjective. 

We also have induced maps $\Phi_m\colon \left(tK\[t\]/t^{m+1}K\[t\]\right)^d\to
\left(tK\[t\]/t^{m+1}K\[t\]\right)^d$. If $T_m$ and $S_m$ are the images of $T$ and $S$,
respectively, in $\left(tK\[t\]/t^{m+1}K\[t\]\right)^d$, and if $m\gg 0$, then $S_m=\Phi_m^{-1}(T_m)$,
and the induced map $S_m\to T_m$ is surjective. 

It follows from Proposition~\ref{lem1_2} that there are $\QQ$-ideals $\frb_+$ and 
$\frb_{-}$ in $\QQ[z_1,\ldots,z_d]$ such that 
$$\max_j\{\ord(\frb(a_1,\ldots,a_d,\alpha_j(a)))\}=
\ord(\frb_+(a_1,\ldots,a_d))-\ord(\frb_-(a_1,\ldots,a_d))$$
for every $a\in \left(K\[t\]\right)^d$. 
 Let $\widetilde{\frb}_+$ and $\widetilde{\frb}_-$
be the $\QQ$-ideals in $\QQ[z_1,\ldots,z_d]$ such that 
$\frb_+(\Phi(\alpha))=\widetilde{\frb}(\alpha)$ and $\frb_-(\Phi(\alpha))
=\widetilde{\frb}_-(\alpha)$. 

\begin{proof}[Proof of Proposition~\ref{subcylinder1}]
By assumption, there is $a\in T$ such that
$\ord({\frb}_+(a))<\infty$ and $\ord({\frb}_-(a))<\infty$.
After replacing $T$ by a subcylinder open in $T$, we may assume that
$\ord(\frb_+(a))=b_+$ and $\ord(\frb_-(a))=b_-$  for all $a\in T$.

For every $j$ and every $\alpha\in S$, we have $\ord(\frb(\Phi(\alpha),\alpha_j))\leq b_+-b_-$,
and the set $S_j$ of those $\alpha$ for which this is an equality is a subcylinder of $S$
that is closed in$S$.
Since $S$ is contained in $\bigcup_jS_j$, 
if $S'$ is an irreducible component of $S$ that dominates $T$, then there is $i$ such that
$S'\subseteq S_i$. 
 Let $S''$ be a subcylinder of $S'$, open in $S'$, and $N_1,\ldots,N_r$ such that
 $\ord(\frq_j(\Phi(\alpha),\alpha_i))=N_j$ for every $\alpha\in S''$, and every $j\leq r$. 
 
 Consider $m\gg 0$, so in particular $m\gg d(b_+-b_-)$ and $m\gg dN_j$ for every $j$. 
 If $S''_m$ is the image of $S''$ in $\left(tK\[t\]/t^{m+1}K\[t\]\right)^d$, then $\Phi_m(S''_m)$
 is constructible, and dense in $T_m$, hence it contains an open subset $T_m^{\circ}$
 of $T_m$. We claim that the inverse image $T^{\circ}$ of $T^{\circ}_m$ in $\left(tK\[t\]\right)^d$
 satisfies the conditions in the proposition.
Indeed, if $a\in T^{\circ}$, then there is $a'\in T$ having the same truncation mod $(t^{m+1})$,
and for which there is a root $\beta'$ of $y^d+\sum_{j=1}^da'_jy^{d-j}$ with
$\ord(\frb(a',\beta'))=b_+-b_-$, and $\ord(\frq_j(a',\beta'))=N_j$ for every $j$.
It follows from Lemma~\ref{lem11} that we can find a root $\beta$ of $y^d+\sum_{j=1}^da_jy^{d-j}$
such that $\ord(\beta-\beta')\geq (m+1)/d$. In this case we also have 
$\ord(\frb(a,\beta))=b_+-b_-$ and $\ord(\frq_j(a,\beta))=N_j$ for every $j$.
\end{proof}

\begin{proof}[Proof of Proposition~\ref{subcylinder2}]
The argument is entirely analogous to the one in the previous proof. Since the set $\widetilde{S}_j$
consisting  of those 
$\alpha\in S$ such that $\ord(\frb(\Phi(\alpha),\alpha_j))\geq L$ is a subcylinder of $S$ that is closed in $S$,
it follows that if $S'$ is an irreducible component of $S$ that dominates $T$, then there is $i$
such that $S'\subseteq\widetilde{S}_i$.
We have a subcylinder $S''$ of $S'$ that is open in $S'$, and $N_j$ for $j\leq r$ such that
$\ord(\frq_j(\Phi(\alpha),\alpha_i))=N_j$ for every $\alpha\in S''$, and every $j\leq r$. 
 Using $S''$ we get a subcylinder $T^{\circ}$ of $T$, open in $T$, following the argument in the proof
 of Proposition~\ref{subcylinder1}.
\end{proof}

\section{The inductive description of singularities of pairs}

In order to motivate our inductive description of log canonical pairs, we start by describing its connection
with Shokurov's Conjecture~\ref{ACC}. Let us fix an algebraically closed field $K$ of characteristic zero.
We first restate Question~\ref{inductive_description1} using the terminology introduced in 
\S 2.

\begin{question}\label{reform}
Given integers $n, d\geq 1$ and a positive rational number $c$, is there a $\QQ$-ideal 
$\frp=\frp_{n,d,c}$ in $K[z_1,\ldots,z_d]$ with the following property:
for every $f=y^d+\sum_{i=1}^da_i(x)y^{d-i}$, with $a_i\in
K\[x_1,\ldots,x_n\]$ of positive order, we have
$$\left({\rm Spec}\,K\[x_1,\ldots,x_n,y\], f^c\right)\,\text{is log canonical if and only if}$$
$$\left({\rm Spec}\,K\[x_1,\ldots,x_n\],\frp(a_1,\ldots,a_d)\right)\,\text{is log canonical}\,?$$
\end{question} 

Before explaining the connection of this question to Conjecture~\ref{ACC}, we make some general considerations. Suppose that $n\geq 1$ is fixed, and consider a nonzero $f\in R=K\[x_1,\ldots,x_n,y\]$,
with $\ord(f)\geq 1$. Set $X={\rm Spec}(R)$. Note that the condition $\lct(f)\geq c$
is independent of the coordinates, or of multiplying $f$ by an invertible element.
If $f(0,\ldots,0,1)$ is nonzero, and its order is $d$, then
by the Weierstrass Preparation Theorem
 $f$ can be written as 
$u\cdot g$, where $u$ is invertible, and $g$ is a Weierstrass polynomial in $y$. Since we are free to replace $f$ by $g$, we may assume that $f$ is a Weierstrass polynomial in $y$, that is,
\begin{equation}\label{eq_Weierstrass}
f=y^d+\sum_{i=1}^da_i(x_1,\ldots,x_n)y^{d-i},
\end{equation}
for some $a_1,\ldots,a_d\in K\[x_1,\ldots,x_n\]$, each having positive order.
Hence, after a linear change of coordinates, we can always assume that $f$ has this form,
and in fact, that $d=\ord(f)$.

\begin{proposition}\label{implication}
Fix $n\geq 1$, and suppose that Question~${\rm \ref{reform}}$ has a positive answer for $n$, and for every
$d$ and $c$. If Conjecture~${\rm \ref{ACC}}$ holds for $n$, then it also holds for $n+1$.
\end{proposition}

\begin{proof}
It is known that the set ${\mathcal HT}_{n+1}$ can also be described as
\begin{equation}\label{set_T}
{\mathcal HT}_{n+1}=\{\lct(f)\mid f\in K\[x_1,\ldots,x_n,y\], f\neq 0, \ord(f)\geq 1\}.
\end{equation}
Furthermore, this set is independent of the algebraically closed field $K$
(for these two facts, see Prop.~3.2 and 3.3 in \cite{dFM}). 
In particular, after possibly extending $K$, we may assume that it is uncountable.

The main idea in the proof is the same as in \emph{loc. cit.}, Prop. 5.2. Suppose that 
we have a strictly  increasing sequence $c_m=\lct(f_m)$, with $f_m\in
K\[x_1,\ldots,x_n,y\]$ nonzero and of positive order. Let $c=\lim_{m\to\infty}c_m$. It follows from
Theorem~1.3 in \emph{loc. cit.} (see also Thm. 29 in \cite{kollar_thr}) that $c\in\QQ$.

Note first that the set of all orders 
$\ord(f_m)$ is bounded. Indeed, each log canonical threshold is positive and since this is an increasing sequence, we can choose a positive $\epsilon$ such that $\epsilon<\lct(f_m)$ for every $m$. On the other hand,
if $f_m\in\frm^N$, where $\frm=(x_1,\ldots,x_n,y)$, then $\lct(f_m)\leq\lct(\frm^N)=\frac{n+1}{N}$. 
We deduce that if $N$ is such that $\frac{n+1}{N}<\epsilon$, then $\ord(f_m)<N$ for every $m$. After possibly replacing our sequence by a subsequence, we may assume that $\ord(f_m)=d$ for every $m$.

It follows as above, by the Weierstrass Preparation Theorem,
that  we may assume
\begin{equation}\label{eq_f_m}
f_m=y^d+\sum_{i=1}^d a_i^{(m)}y^{d-i},
\end{equation}
with $a_i^{(m)}\in K\[x_1,\ldots,x_n\]$ of positive order.
For every positive rational number $\tau$, let 
$${\mathcal A}_{n,d}(\tau):=\left\{(a_1,\ldots,a_d)\in \left(K\[x_1,\ldots,x_n]\right)^d\mid \ord(a_i)\geq 1\,\text{for all}\,i,\,
\lct\left(y^d+\sum_{i=1}^da_iy^{d-i}\right)\geq \tau\right\}.$$

In order to get a contradiction, it is enough to show that every such set is a cylinder, that is, it is the inverse image of a constructible set by some projection $\left(K\[x_1,\ldots,x_n\]\right)^d\to \left(K\[x_1,\ldots,x_n\]/(x_1,\ldots,x_n)^s\right)^d$. Indeed, then we have the sequence of cylinders
$${\mathcal A}_{n,d}(c_1)\supset {\mathcal A}_{n,d}(c_2)\supset \ldots\supset
\bigcap_m{\mathcal A}_{n,d}(c_m)={\mathcal A}_{n,d}(c).$$
By assumption, all these inclusions are strict, and since $K$ is assumed uncountable, this contradicts
a basic property of cylinders: a decreasing sequence of cylinders with empty intersection is eventually empty (the proof in our context follows verbatim the proof of Lemma~1.2 in \cite{ELM}). 

Suppose now that that Question~\ref{reform} has a positive answer for $n$, $d$ and $\tau$,
with $\frp=(P_1,\ldots,P_r)^{q}$.
We can then rewrite 
$${\mathcal A}_{n,d}(\tau)=\left\{a=(a_1,\ldots,a_d)\in \left(K\[x_1,\ldots,x_n\]\right)^d\mid
\ord(a_i)\geq 1\,\text{for all}\,i, \lct(P_1(a),\ldots,P_r(a))\geq q\right\},$$
and we need to show that if Conjecture~\ref{ACC} holds in dimension $n$, then this set is a cylinder.
On the other hand, the map 
$$(P_1,\ldots,P_r)\colon \left((x_1,\ldots,x_n)K\[x_1,\ldots,x_n\]\right)^d\to 
\left(K\[x_1,\ldots,x_n\]\right)^r$$
pulls-back cylinders to cylinders, hence it is enough to show that the set
$${\mathcal B}_{n,r}(q)=\{(g_1,\ldots,g_r)\in \left(K\[x_1,\ldots,x_n\]\right)^r\mid
\lct(g_1,\ldots,g_r)\geq q\}$$
is a cylinder in $\left(K\[x_1,\ldots,x_n\]\right)^r$. 

Since Conjecture~\ref{ACC} holds in dimension $n$, the set
$${\mathcal T}_n:=\left\{\lct(J)\mid J\subset K\[x_1,\ldots,x_n\], \ord(J)\geq 1, J\neq (0)\right\}$$
also satisfies ACC. Indeed, every element in this set can be written as $n\cdot\lct(g)$, for some
$g\in K\[x_1,\ldots,x_n\]$ (for this fact, see Prop. 3.5 and 3.3 in \cite{dFM}). Therefore we can find
$\delta>0$ such that ${\mathcal T}_n$ contains no real number in the interval $(q-\delta,q)$. 
On the other hand, it follows from  Corollary~2.11 in \emph{loc. cit.} that if $M$ is a positive integer such that $\frac{n}{M}<\delta/2$, and if 
$J_1$ and $J_2$ are ideals in $K\[x_1,\ldots,x_n\]$ that induce the same ideal in 
$K\[x_1,\ldots,x_n\]/(x_1,\ldots,x_n)^M$, then $\lct(J_1)-\lct(J_2)|\leq\frac{n}{M}<\delta/2$.
In particular, we see that $\lct(J_1)\geq q$ if and only if $\lct(J_2)\geq 2$. Therefore 
the set ${\mathcal B}_{n,r}(q)$ is the inverse image of a set in 
$\left(K\[x_1,\ldots,x_n\]/(x_1,\ldots,x_n)^M\right)^r$. The fact that this set is open is then a consequence of the semicontinuity of log canonical thresholds (see, for example, Thm.4.9 in \cite{Mus}). We conclude
that ${\mathcal B}_{n,r}(q)$, hence also ${\mathcal A}_{n,d}(\tau)$, is a cylinder. This completes
the proof of the proposition.
\end{proof}

\bigskip

We now state our main result, but 
let us first set some notation. Suppose that $c$ is a positive rational number. We are interested in
whether $\lct(f)\geq c$ for a given formal power series $f$ with $1\leq \ord(f)\leq d$. Note that for every such $f$ we have $\frac{1}{d}\leq\lct(f)\leq 1$ (the first inequality is well-known in the finite type case,
see Lem. 8.10 in \cite{kollar}, and in our case can be deduced via Prop. 2.5 in \cite{dFM}).
Therefore we may assume that $\frac{1}{d}<c\leq 1$. Let $p\in\{1,\ldots,d-1\}$ be such that
\begin{equation}\label{ineq_c}
\frac{1}{d-p+1}<c\leq\frac{1}{d-p}.
\end{equation}
We apply Corollary~\ref{cor_section2} with $k=p$, $c_1=1-(d-p)c$, and $c_2=(d-p+1)c-1$
to get the $\QQ$-ideals $\frp_+$ and $\frp_{-}$ in $\QQ[z_1,\ldots,z_d]$.
Note that our assumptions imply $c_1\geq 0$ and $c_2>0$.

\begin{theorem}\label{main}
Fix an integer $d\geq 1$, an integer $p$ with $1\leq p\leq d-1$, and a rational number $c$ satisfying 
${\rm (6)}$. For every ${\rm (}$algebraically closed, characteristic zero${\rm )}$ field $K$, and every
$a_1,\ldots,a_d\in  K\[x_1,\ldots,x_n\]$ of positive order, if $f=y^d+\sum_{i=1}^da_iy^{d-i}$, then
$\lct(f)\geq c$ if and only if $({\rm Spec}\,K\[x_1,\ldots,x_n\],
\frp_+(a_1,\ldots,a_d)\frp_{-}
(a_1,\ldots,a_d)^{-1})$ is log canonical. 
\end{theorem}

\begin{remark}\label{inclusion_ideals}
It follows from our definition and Corollary~\ref{cor_section2} i) that for every $K$ we have 
$$\frp_+ K[z_1,\ldots,z_d]\subseteq \left(\frp_- K[z_1,\ldots,z_d]\right)^{\frac{d}{d-1}}.$$ 
\end{remark}

\begin{remark}
Note that by definition, the condition that the pair 
$$({\rm Spec}\,K\[x_1,\ldots,x_n\],
\frp_+(a_1,\ldots,a_d)\frp_{-}
(a_1,\ldots,a_d)^{-1})$$ in the theorem be log canonical implies that the two $\QQ$-ideals
$\frp_+(a_1,\ldots,a_d)$ and $\frp_{-}(a_1,\ldots,a_d)$ are nonzero.
\end{remark}

\begin{example}\label{degree_three}
Let us now describe explicitly the case $d=3$. After a  change of coordinates
of the form $y'=y-g(x)$, we may write
$f=y^3+a(x)y+b(x)$, with $a$, $b\in K\[x_1,\ldots,x_n\]$ of positive order. Recall
 that we denote by $\Delta_f$ the discriminant $4a(x)^3+27b(x)^3$. 
It follows from Theorem~\ref{main} and Example~\ref{example1}
 that $\lct(f)\geq c$ if and only if 
$$\left\{
\begin{array}{cl}
\left({\rm Spec}\,K\[x_1,\ldots,x_n\], \left(a(x)^3, b(x)^2\right)^{\frac{3c-1}{6}}\right)\,
\text{is log canonical}\,, & {\rm when}\,\frac{1}{3}<c\leq\frac{1}{2}; \\[2mm]
\left({\rm Spec}\,K\[x_1,\ldots,x_n\], (\Delta_f)^{c-\frac{1}{2}}\cdot \left(a(x)^3, b(x)^2\right)^{\frac{2-3c}{6}}\right)
\,\text{is log canonical}\,, & {\rm when}\,\frac{1}{2}<c\leq 1.
\end{array}\right.
$$
Note that the pair that appears in our condition is not effective if $c\in\left(\frac{2}{3},1\right]$.
\end{example}

\begin{example}\label{example3}
A simple case is when $\frac{1}{d}<c\leq\frac{1}{d-1}$. 
After a change of variable, we may assume that $a_1=0$. 
In this case Theorem~\ref{main}
and Example~\ref{example0} imply that if $f=y^d+\sum_{i=2}^da_iy^{d-i}$, 
with all $a_i\in K\[x_1,\ldots,x_n\]$ of positive order, then 
$\lct(f)\geq c$ if and only if 
$$({\rm Spec}\,K\[x_1,\ldots,x_n\], \frp(a_2(x),\ldots,a_d(x))^{cd-1})$$
is log canonical, where $\frp=\sum_{i=2}^d(z_i)^{1/i}$.
\end{example}

\begin{remark}
Recall that our goal is to understand the set
$$\left\{(a_1,\ldots,a_d)\in \left(K\[x_1,\ldots,x_n\]\right)^d\mid\ord(a_i)\geq 1\,\text{for all}\,i,\,
\lct\left(y^d+\sum_{i=1}^da_iy^{d-i}\right)
\geq \tau\right\},$$ where $\tau$ is a fixed positive rational number.
Since $\lct(f^m)=\lct(f)/m$,  after replacing $d$ and $\tau$ by $md$ and $\tau/m$, for a suitable $m$, we see that
it is enough to describe such a set when
 $\tau$ is of the form $\frac{1}{q}$, with $q\in\{1,\ldots,d-1\}$.
In this case, in Theorem~\ref{main} we have $c_1=0$, and the condition 
for $\frp_+$ and $\frp_-$ becomes substantially simpler. 
\end{remark}

We now show that the statement of Theorem~\ref{main} is strong enough to give Conjecture~\ref{ACC}
in dimension two (see \cite{Sho}, \cite{PS}, and \cite{FJ} for other proofs of this case of the conjecture).

\begin{corollary}\label{dimension_two}
Conjecture~${\rm \ref{ACC}}$ holds for $n=2$.
\end{corollary}

\begin{proof}
Arguing as in the proof of Proposition~\ref{implication}, we see that  it is enough to show the following:
given $d\geq 1$, a positive rational number $\tau$, and an uncountable algebraically closed
field $K$ of characteristic zero, the set
$${\mathcal A}_{2,d}(\tau):=\left\{(a_1,\ldots,a_d)\in \left(K\[x\]\right)^d\mid \ord(a_i)\geq 1\,\text{for all}\,i,\,\lct\left(y^d+\sum_{i=1}^da_iy^{d-i}\right)\geq \tau\right\}$$
is a cylinder in $\left(K\[x\]\right)^d$. Let us apply Theorem~\ref{main} 
for $d$ and $\tau$ to get the $\QQ$-ideals $\frp_+=(P_1,\ldots,P_r)^{p}$ and $\frp_{-}=
(Q_1,\ldots,Q_s)^q$ in $\QQ[z_1,\ldots,z_d]$. By Remark~\ref{inclusion_ideals} we have
\begin{equation}\label{eq_dim_two}
\frp_+ K[z_1,\ldots,z_d]\leq \left(\frp_- K[z_1,\ldots,z_d]\right)^{\frac{d}{d-1}}.
\end{equation}

Note that if $\fra_1,\fra_2$ are $\QQ$-ideals in $K\[x\]$, we have 
$({\rm Spec}\,K[x],\fra_1\fra_2^{-1})$ log canonical if and only if $\ord(\fra_1)-\ord(\fra_2)\leq 1$.
Therefore we can describe
${\mathcal A}_{2,d}(\tau)$ as the set of those $a=(a_1,\ldots,a_d)\in \left(K\[x\]\right)^d$ with $\ord(a_i)\geq 1$ for all $i$, and such that 
\begin{equation}\label{eq2_dim_two}
p\cdot \ord(P_1(a),\ldots,P_r(a))-q\cdot\ord(Q_1(a),\ldots,Q_s(a))\leq 1.
\end{equation}
On the other hand, (\ref{eq_dim_two}) gives
$$p\cdot\ord(P_1(a),\ldots,P_r(a))\geq\frac{dq}{d-1}\cdot \ord(Q_1(a),\ldots,Q_s(a))$$
for every $a\in \left(K\[x\]\right)^d$.
This implies that in order to have (\ref{eq2_dim_two}), we need 
$\ord(Q_1(a),\ldots,Q_s(a))\leq \frac{d-1}{q}$, and therefore also 
$\ord(P_1(a),\ldots,P_r(a))\leq \frac{d}{p}$. In particular, there are only finitely many possibilities
for  $\ord(P_1(a),\ldots,P_r(a))$ and $\ord(Q_1(a),\ldots,Q_s(a))$ such that
(\ref{eq2_dim_two}) holds. Since for every pair $(i,j)$ the set
$$\{a=(a_1,\ldots,a_d)\in (K\[x\])^d\mid \ord(P_1(a),\ldots,P_r(a))=i,\,
\ord(Q_1(a),\ldots,Q_s(a))=j\}$$
is clearly a cylinder, it follows that ${\mathcal A}_{2,d}(\tau)$ is a cylinder.
\end{proof}

\bigskip

Before giving the proof of Theorem~\ref{main}, we make some preparations.
In light of Theorem~\ref{description}, in order to decide whether $\lct(f)\geq c$, we need to 
estimate the codimension of the set of those $(u_1,\ldots,u_n,w)\in (t K\[t\])^{n+1}$
with $\ord(f(u_1,\ldots,u_n,w))=m$. The key idea in our approach is to first plug in
$n$-uples $u=(u_1,\ldots,u_n)\in (tK\[t\])^n$ to get $h(y):=f(u,y)\in K\[t\][y]$, and then analyze the condition
on $w\in tK\[t\]$ to get $\ord(h(w))=m$. Note that $h$ is a polynomial with formal power 
series coefficients to which we can apply the considerations in the previous section.

We now fix $f$ given by (\ref{eq_Weierstrass}), with the $a_i$ in $K\[x_1,\ldots,x_n\]$
of positive order. At this point let us fix also $u=(u_1,\ldots,u_n)
\in (tK\[t])^n$, and let $\alpha_1,\ldots,\alpha_d$ be the roots of $h=f(u,y)$, hence
$$h=\prod_{i=1}^d(y-\alpha_i).$$
Using the notation in Corollary~\ref{cor3}, we assume that the order of  $P_i(a_1(u),\ldots,a_d(u))$
is divisible by $m$ for $i\leq r$ (due to the second assertion in Theorem~\ref{description},
this assumption will turn out to be harmless). Therefore $\alpha_1,\ldots,\alpha_d\in K\[t\]$.
In fact, they lie in $tK\[t\]$, since they have positive order by Lemma~\ref{lem1}. 

\begin{lemma}\label{lem3_1}
With the above notation, the following hold:
\begin{enumerate}
\item[i)] For every $i$ with $1\leq i\leq d$, and every $w\in tK\[t\]$, we have
\begin{equation}\label{eq_lem3_1}
\ord(f(u,w))\geq\sum_{j=1}^d\min\{\ord(w-\alpha_i), \ord(\alpha_i-\alpha_j)\}.
\end{equation}
\item[ii)] If $i$ is such that $\ord(w-\alpha_i)=\max_j\{\ord(w-\alpha_j)\}$, then the
inequality in ${\rm (\ref{eq_lem3_1})}$ is an equality.
\item[iii)] For every $i$ with $1\leq i\leq d$, and every $q\geq 1$, the set of those $w\in tK\[t\]$ with 
$q=\ord(w-\alpha_i)=\max_j\{\ord(w-\alpha_j)\}$ is a locally closed cylinder in
$K\[t\]$ of codimension $q$.
\end{enumerate}
\end{lemma}

\begin{proof}
It is clear that $\ord(f(u,w))=\sum_{j=1}^d\ord(w-\alpha_j)$. Since 
\begin{equation}\label{eq2_lem3_1}
\ord(w-\alpha_j)\geq\min\{\ord(w-\alpha_i),\ord(\alpha_i-\alpha_j)\},
\end{equation}
 we get i). Moreover, if $i$ is such that $\ord(w-\alpha_i)\geq\ord(w-\alpha_j)$
 for every $j$, then in (\ref{eq2_lem3_1}) we have equality.  This gives ii).
 For iii), note that if $\alpha_j=\sum_{m\geq 1}c_{j,m}t^m$ and if we write $w=\sum_{m\geq 1}c_mt^m$, then  the condition $\ord(w-\alpha_i)=q\geq\ord(w-\alpha_j)$
 for every $j$ is equivalent to 
 $$c_m=c_{i,m}\,\text{for}\,m\leq q-1,\,\text{and}\, c_q\neq c_{j,q}\,\text{for those}\,j\,\text{with}
 \ord(\alpha_i-\alpha_j)\geq q.$$
 The set of such $w$ is a locally closed cylinder in $K\[t\]$ of codimension $q$.
\end{proof}

\bigskip

We can now give the proof of our main result.

\begin{proof}[Proof of Theorem~\ref{main}]
Fix $K$ and $a_1,\ldots,a_d\in K\[x_1,\ldots,x_n\]$ of positive order.
 Let $f=y^d+\sum_{i=1}^da_iy^{d-i}$. For every $u\in \left(tK\[t\]\right)^n$, we denote by
 $\alpha_1(u),\ldots,\alpha_d(u)$ the roots of $f(u,y)$.
 Since a pair is log canonical over $K$ if and only if it is log canonical over a field extension of $K$,
 we may assume that $K$ is uncountable. 

We first assume that the pair $({\rm Spec}\,K\[x_1,\ldots,x_n\],
\frp_+(a_1,\ldots,a_d)\frp_{-}
(a_1,\ldots,a_d)^{-1})$ is log canonical, and prove that $\lct(f)\geq c$. 
 By Theorem~\ref{description}, it is enough to show that for every irreducible locally closed
cylinder $C\subseteq (tK\[t\])^{n+1}$, we have $\codim(C)\geq c\cdot\ord_C(f)$. Of course, in order to check this inequality
we can always replace $C$ by a suitable subcylinder that is open in $C$. 

It follows from the statement of 
Theorem~\ref{description} that given finitely many nonzero power series in 
$K\[x_1,\ldots,x_{n+1}\]$, 
after possibly replacing $C$ by a subcylinder open in $C$, 
we may assume that for every $\gamma\in C$,
the order of $\gamma$ along each of these power series is constant, and divisible enough. 
In particular, we may assume that for every $(u,w)\in (tK\[t\])^n\times tK\[t\]$ that lies in $C$, we have 
$\ord(f(u,w))=m=\ord_C(f)$, and $\ord(\frp_+(a_1(u),\ldots,a_d(u)))=\beta_+$
and $\ord(\frp_-(a_1(u),\ldots,a_d(u)))=\beta_-$.
By a similar argument, using also Corollary~\ref{cor3}, we may assume that for every $(u,w)\in C$,
  all $\alpha_i(u)$ lie in $K\[t\]$. Furthermore, applying Corollary~\ref{cor2} for $k=d$ and 
  $k=(d-1)$, we may also assume that for all $(u,w)\in C$, we have $\max_j\left\{\ord(\alpha_j(u)-w)
  \right\}=q$,
  for some $q\in\ZZ_{\geq 0}$ (note that $q$ is finite: otherwise $f(u,w)=0$ for every $(u,w)\in C$,
  which would force $f$ to be zero). 
  
  Let now $N$ be large enough, so that $C$ is the inverse image of $C_N\subseteq
(tK\[t\]/t^NK\[t\])^{n+1}$. Let $C'_N\subseteq (tK\[t\]/t^NK\[t\])^{n}$ be the image of $C_N$
by the map induced by the projection onto the first $n$ components.  
Since we may assume that $N>q$, it follows from Lemma~\ref{lem3_1} iii) that 
for every ${u}\in C'_N$, the fiber of $C_N$ over ${u}$ has dimension
$\leq N-q$. Therefore we have 
$$\codim(C)=N(n+1)-\dim(C_N)\geq N(n+1)-(N-q)-\dim(C'_N).$$
It follows that in order to guarantee $\codim(C)\geq cm$ it would be enough to have 
\begin{equation}\label{formula13}
\codim(C')\geq cm-q, 
\end{equation} 
where $C'\subseteq (tK\[t\])^n$ is the inverse image of $C_N'$. 

\noindent{\bf Claim}. If $u\in C'$, then 
\begin{equation}\label{formula14}
\max_{i=1}^d\left\{c_1\cdot \min_{j_1,\ldots,j_{p-1}}\ord\left(\prod_{\ell=1}^{p-1}(\alpha_{j_{\ell}}(u)-\alpha_i(u))\right)+c_2\cdot \min_{j_1,\ldots,j_{p}}\ord\left(\prod_{\ell=1}^{p}(\alpha_{j_{\ell}}(u)-
\alpha_i(u))\right)
\right\}\geq cm-q
\end{equation}
(recall that $c_1=1-(d-p)c$ and $c_2=(d-p+1)c-1$). Indeed, note first that if the left-hand side of
(\ref{formula14}) is infinite, then there is nothing to prove. Otherwise, this is equal to
$\ord\left(\frp_+(a_1(u),\ldots,a_d(u)\right))-\ord\left(\frp_-(a_1(u),\ldots,a_d(u))\right)$,
and both these orders are finite. We see that if $N$ is larger than $\beta_+, \beta_-, cm-q$, then in order to check (\ref{formula14})
we may replace $u$ by some $u'$ that has the same image modulo $t^N$, hence we may assume that
for some $w$ we have $(u,w)\in C$. 

Let $i$ be such that $\ord(w-\alpha_i(u))=q$.  It follows from Lemma~\ref{lem3_1} that 
  \begin{equation}\label{formula11}
  m=\ord(f(u,w))=\sum_{j=1}^d\min\{\ord(\alpha_i(u)-\alpha_j(u)),q\}.
  \end{equation}
Let us order the $\ord(\alpha_i(u)-\alpha_j(u))$, for $1\leq j\leq d$, as $b_1\leq\ldots\leq b_d=\infty$.
Suppose that $r$ is such that $b_r<\infty$, and $b_{r+1}=\infty$ (hence $0\leq r\leq d-1$,
with the convention that $b_0=1$). Let $k\in\{0,\ldots,r\}$ be largest with $b_k\leq q$. 
We deduce  from (\ref{formula11}) that 
$m=b_1+\ldots+b_k+(d-k)q$.
Note that $\beta_+-\beta_-\geq c_1(b_1+\cdots+b_{p-1})+c_2(b_1+\cdots+b_p)$, and since
$c_2>0$, we have $b_p<\infty$.
Moreover, in order to show (\ref{formula14}) it is enough to show 
$$c_1(b_1+\cdots+b_{p-1})+c_2(b_1+\cdots+b_{p})\geq cm-q,$$
which is equivalent to 
\begin{equation}\label{formula15}
c(b_1+\cdots+b_{p-1})+\left(c(d-p)-1\right)b_{p}\geq c(b_1+\cdots+b_k)+\left(c(d-k)-1\right)q.
\end{equation}

By our choice of $p$, the right-hand side of (\ref{formula15}) is 
weakly increasing in $q$ for $k\leq p-1$, hence for 
$q\in [b_k,b_{k+1})$ it is $\leq c(b_1+\cdots+b_k)+\left((d-k)c-1\right)b_{k+1}$.
This expression, in turn, is weakly increasing in $k$ when $k\leq p-1$: the inequality
$$c(b_1+\cdots+b_k)+\left(c(d-k)-1\right)b_{k+1}\leq c(b_1+\cdots+b_{k+1})+\left(c(d-k-1)-1\right)
b_{k+2}$$
follows from $b_{k+1}\leq b_{k+2}$ and $c(d-k-1)-1\geq 0$, which holds
for $k\leq p-2$.
 Hence we get (\ref{formula15})
when $k\leq p-1$. A similar argument gives (\ref{formula15}) also for $k\geq p$:
in this case, the right-hand side of (\ref{formula15}) is weakly decreasing in $q$, hence it is
$\leq c(b_1+\ldots+b_{k-1})+(c(d-k+1)-1)b_k$. On the other hand, for $k\geq p$ this expression
is weakly decreasing in $k$, hence attains its maximum for $k=p$.
  This completes the proof of our claim. 

 By definition, the left-hand side of (\ref{formula14}) is equal to
$$\ord_{C'}(\frp_+(a_1(u),\ldots,a_d(u)))-\ord_{C'}(\frp_{-}(a_1(u),\ldots,a_d(u)).$$ 
Since the pair
$({\rm Spec}\,K\[x_1,\ldots,x_n\],\frp_+(a_1,\ldots,a_d)\frp_{-}
(a_1,\ldots,a_d)^{-1})$ is log canonical, Theorem~\ref{description} and the above claim imply that
$\codim(C')\geq cm-q$, and as we have seen, we get $\codim(C)\geq c\cdot\ord_C(f)$.
Since this holds for every $C$ as above, we conclude that $\lct(f)\geq c$. 

We now prove the converse. Suppose that $\lct(f)\geq c$. Arguing by contradiction, 
let us assume that the pair $({\rm Spec}\,K\[x_1,\ldots,x_n\],
\frp_+(a_1,\ldots,a_d)\frp_{-}
(a_1,\ldots,a_d)^{-1})$ is not log canonical. 
We first consider the case when $\frp_+(a_1,\ldots,a_n)$
(hence also $\frp_{-}(a_1,\ldots,a_n)$) is nonzero. By definition and Theorem~\ref{description}, we can find an irreducible locally closed cylinder 
$D'\subseteq \left(t K\[t\]\right)^n$ such that 
for every  $u\in D'$ we have $\ord(\frp_{+}(a_1(u),\ldots,a_d(u)))=m_1$ and
$\ord(\frp_{-}(a_1(u),\ldots,a_d(u)))=m_2$, and $\codim(D')<m_1-m_2$. Furthermore, given 
finitely many power series in $K\[x_1,\dots,x_{n}\]$, we may assume that the order of each of them
is constant along the $u\in D'$, and is divisible enough. In particular, if we denote as above by
$\alpha_1(u),\ldots,\alpha_d(u)$ the roots of $f(u,y)$, we may assume by Corollary~\ref{cor3}
that all $\alpha_i(u)\in K\[t\]$ for $u\in C'$ (in which case they lie in $tK\[t\]$ by
Lemma~\ref{lem1}). 

As before, given $i\leq d$, for every $u\in D'$ we order the $\ord(\alpha_j(u)-\alpha_i(u))$
with $1\leq j\leq d$ 
as $b_1(u)\leq\ldots\leq b_d(u)=\infty$. By Corollary~\ref{cor100} we have $\QQ$-ideals
$\frp_+^{(k)}$ and $\frp_-^{(k)}$ such that 
$$b_k(u)=\ord(\frp_+^{(k)}(a_1(u),\ldots,a_d(u),\alpha_i(u)))-
\ord(\frp^{(k)}_-(a_1(u),\ldots,a_d(u),\alpha_i(u))).$$
We have seen in the proof of Corollary~\ref{cor_section2} that there 
is a $\QQ$-ideal $\frp$ in $\QQ[z_1,\ldots,z_{d+1}]$ such that 
$$\ord(\frp(a(u),\alpha_i(u)))=c_1\cdot \min_{j_1,\ldots,j_{p-1}}\ord\left(\prod_{\ell=1}^{p-1}(\alpha_{j_{\ell}}(u)-\alpha_i(u))\right)+c_2\cdot \min_{j_1,\ldots,j_{p}}\ord\left(\prod_{\ell=1}^{p}(\alpha_{j_{\ell}}-\alpha_i)\right).$$
Applying Proposition~\ref{subcylinder1} to $\frp$, and to all $\frp^{(k)}_+$, $\frp_-^{(k)}$, we see 
that after possibly replacing $D'$ by a subcylinder open in $D'$, we may assume that there are 
$N_1,\ldots,N_d$
such that for every $u\in D'$ we can find an $i$ with 
\begin{equation}\label{eq_thm_main1}
\ord(\frp(a(u),\alpha_i(u)))=\max_j\{\ord(\frp(a(u),\alpha_i(u)))\},\,\,
b_k(u)=N_k
\end{equation}
for $1\leq k\leq d$ (where, of course, $b_k(u)$ is computed with respect to this $i$).

On the other hand, note that for such $i$ we have
\begin{equation}\label{eq_thm_main2}
c\cdot\sum_{k=1}^{p-1}b_k(u)+(c(d-p+1)-1)b_p(u)=c_1\cdot \sum_{k=1}^{p-1}b_k(u)+c_2\cdot \sum_{k=1}^{p}b_k(u)=m_1-m_2.
\end{equation}
In particular, this implies $N_1,\ldots,N_p<\infty$.
Let $m=N_1+\ldots+N_{p-1}+(d-p+1)N_p$, and consider the locally closed cylinder 
$$D:=\{(u,w)\in D'\times tK\[t\]\mid \ord(f(u,w))\geq m\}.$$
Let $q=N_p$, and consider the fiber of the projection $D\to D'$
over some $u\in D'$. We have seen that there is $\alpha_i(u)$ such that
(\ref{eq_thm_main1}) holds. By Lemma~\ref{lem3_1}, the set
$$D(u):=\{w\in tK\[t\]\mid \ord(w-\alpha_i(u))=q\}$$ is a locally closed cylinder 
(coming from a level independent on $u$)
such that $D(u)$ is contained in the fiber of $D$ over $u\in D'$, and
$\codim(D(u))=q$. We conclude that there is an irreducible component 
of $D$ having codimension 
$\leq q+\codim(D')<N_p+m_1-m_2=cm$, where we used the formula (\ref{eq_thm_main2}).
This contradicts the fact that $\lct(f)\geq c$.

Suppose now that $\frb_+(a_1,\ldots,a_d)=0$. Let us consider an arbitrary irreducible locally closed
cylinder $D'\subseteq \left(tK\[t\]\right)^n$, such that for every $u\in D'$, all roots 
$\alpha_i(u)$ are in $K\[t\]$. Let $\frp$, $\frp_{+}^{(k)}$, and $\frp_-^{(k)}$ be as above. 
By Proposition~\ref{subcylinder2}, given $L>0$, after replacing $D'$ by a subcylinder open in $D'$, we may assume that there are $N_1,\ldots,N_d$ such that for every $u\in D'$ we can find $i$ with
$$\ord(\frp(a(u),\alpha_i(u)))\geq L,\,\,b_k(u)=N_k,$$
for all $k$.
Note that this implies, in particular, $c\cdot\sum_{k=1}^{p-1}N_k+(c(d-p+1)-1)N_p\geq L$.

Let $q\geq 1$ be fixed, and $m=\sum_{k=1}^d\min\{q,N_k\}$. Consider
$$D:=\{(u,w)\in D'\times tK\[t\]\mid \ord(f(u,w))\geq m\}.$$
Arguing as before, we see that $\lct(f)\geq c$ implies
$\codim(D')+q\geq c\cdot m$. By letting $q$ go to infinity, we deduce
$1\geq c\cdot\#\{j\mid N_j=\infty\}$. Therefore at most 
$(d-p)$ of the $N_j$ are infinite, which implies $N_p<\infty$. 
By taking $q=N_p$, we conclude as above that
$$\codim(D')\geq c\cdot\sum_{k=1}^{p-1}N_k+(c(d-p+1)-1)N_p\geq L.$$
Since this holds for every $D'$ and every $L$, we get a contradiction. This completes the 
proof of the theorem.
\end{proof}

\providecommand{\bysame}{\leavevmode \hbox \o3em
{\hrulefill}\thinspace}


\begin{thebibliography}{ELM}


\bibitem[Ale]{alexeev}
V.~A.~Alexeev, Two two-dimensional terminations, Duke Math. J.
\textbf{69} (1993), 527--545.

\bibitem[Bir]{birkar}
C.~Birkar, Ascending chain condition for log canonical thresholds
and termination of log flips, Duke Math. J. \textbf{136} (2007),
173--180.

\bibitem[dFM]{dFM}
T.~de Fernex and M.~Musta\c{t}\u{a}, Limits of log canonical thresholds,
 Ann. Sci. \'{E}cole Norm. Sup., to appear.

\bibitem[ELM]{ELM}
L.~Ein, R.~Lazarsfeld and M.~Musta\c{t}\u{a}, Contact loci in arc
spaces, Compos. Math. \textbf{140} (2004), 1229--1244.

\bibitem[EM]{EM}
L.~Ein and M.~Musta\c{t}\u{a}, Invariants of singularities of pairs,
in \emph{International Congress of Mathematicians. Vol. II}, 583--602, 
Eur. Math. Soc., Z\"{u}rich, 2006. 

\bibitem[FJ]{FJ}
C.~Favre and M.~Jonsson,
Valuations and multiplier ideals. 
J. Amer. Math. Soc. 18 (2005), 655--684.

\bibitem[Kaw]{Kawakita}
M.~Kawakita,  On a comparison of minimal log discrepancies in terms of motivic integration,
J. Reine Angew. Math. \textbf{620} (2008), 55--65.

\bibitem[Kob]{koblitz}
N.~Koblitz, $P$-adic numbers, $p$-adic analysis, and zeta-functions.
Second edition. Graduate Texts in Mathematics \textbf{58},
Springer-Verlag, New York, 1984.

\bibitem[Kol1]{kollar_thr}
J.~Koll\'{a}r, Which powers of holomorphic functions are integrable?,
preprint, arXiv: 0805.0756.

\bibitem[Kol2]{kollar}
J.~Koll\'{a}r, Singularities of pairs, in \emph{Algebraic geometry,
Santa Cruz 1995}, 221--286, volume \textbf{62} of Proc. Symp. Pure
Math. Amer. Math. Soc. 1997.

\bibitem[Kuw]{Kuwata}
T.~Kuwata,  On log canonical thresholds of reducible plane curves, Amer. J. Math. \textbf{121} (1999), 701--721.

\bibitem[Laz]{positivity}
 R.~Lazarsfeld,
\emph{Positivity in algebraic geometry} II, Ergebnisse der
Mathematik und ihrer Grenzgebiete, vol. \textbf{49},
Springer-Verlag, Berlin, 2004.

\bibitem[Mat]{Matsumura}
H.~Matsumura,  \emph{Commutative algebra}, Second edition,
 Mathematics Lecture Note Series \textbf{ 56},
 Benjamin/Cummings Publishing Co., Inc., Reading, Mass., 1980.



\bibitem[Mus]{Mus}
M. Musta\c t\u a, Singularities of pairs via jet schemes, J. Amer.
Math. Soc \textbf{15} (2002), 599--615.

\bibitem[PS]{PS}
D.-H. Phong and J. Sturm, On a conjecture of Demailly and Koll\'{a}r,
Asian J. Math. \textbf{4} (2000), 221--226. 

\bibitem[Sho]{Sho}
V.~V.~Shokurov,
 Three-dimensional log perestroikas. With an appendix
in English by Yujiro Kawamata. {Izv. Ross. Akad. Nauk Ser. Mat.}
\textbf{56} (1992), 105--203; translation in {Russian Acad. Sci.
Izv. Math.} \textbf{40} (1993), 95--202.

\bibitem[Tem]{Temkin}
M.~Temkin,  Desingularization of quasi-excellent schemes in
characteristic zero, Adv. Math. \textbf{219} (2008), 488--522.

\end{thebibliography}

\end{document}